\documentclass[12pt]{amsart}
\usepackage{amssymb}
\usepackage{amscd}
\usepackage{amsmath}

\newcommand{\gog}{{\mathfrak g}}
\newcommand{\gb}{{\mathfrak b}}
\newcommand{\gh}{{\mathfrak h}}
\newcommand{\gl}{{\mathfrak l}}
\newcommand{\gm}{{\mathfrak m}}
\newcommand{\gn}{{\mathfrak n}}

\newcommand{\gs}{{\mathfrak s}}
\newcommand{\bB}{{\bf B}}
\newcommand{\bG}{{\bf G}}
\newcommand{\bM}{{\bf M}}
\newcommand{\bP}{{\bf P}}
\newcommand{\bS}{{\bf S}}
\newcommand{\bT}{{\bf T}}

\newcommand{\al}{\alpha}
\newcommand{\be}{\beta}
\newcommand{\de}{\delta}

\newcommand{\si}{\sigma}
\newcommand{\lam}{\lambda}

\newcommand{\cal}{\mathcal} 
\newcommand{\Cscr}{{\cal C}}
\newcommand{\Oscr}{{\cal O}}
\newcommand{\Vscr}{{\cal V}}
\newcommand{\Wscr}{{\cal W}}
\newcommand{\Dscr}{{\cal D}}
\newcommand{\Iscr}{{\cal I}}
\newcommand{\Yscr}{{\cal Y}}
\newcommand{\Lie}{{\rm Lie\,}}
\newcommand{\sh}{{\rm sh\,}}
\newcommand{\codim}{{\rm codim\,}}
\newcommand{\snext}{{\rm snext}}
\newcommand{\sprev}{{\rm sprev}}
\newcommand{\prev}{{\rm prev}}
\newcommand{\next}{{\rm next}}
\newcommand{\QED}{\par \hspace{15cm}$\blacksquare$ \par}
\newcommand{\pr}{^{\prime}}
\newcommand{\prpr}{^{\prime\prime}}
\newcommand{\st}{\subset}
\newcommand{\Co}{{\mathbb C}} 
\newcommand{\Na}{{\mathbb N}}  
\newcommand{\boT}{{\mathbb T}}
\newcommand{\boS}{{\mathbb S}}

\newcommand{\Pf}{\noindent{\bf Proof.}\par\noindent}
\newcommand{\sr}{\scriptscriptstyle}
\newcommand{\vb}{\vrule height 14pt depth 7pt} 
\newcommand{\ts}{\tabskip 4pt}	
\newcommand{\vsa}{\noalign{\vskip-7pt}}
\newcommand{\ssa}{\noalign{\vskip -1pt}}

\newcommand{\ov}{\overline}
\newcommand{\parno}{\par\noindent}
\newcommand{\dor}{\stackrel{\rm D}{\leq}}     
\newcommand{\dg}{\stackrel{\rm D}{\geq}}

\newcommand{\dgs}{\stackrel {\rm D}{>}}
\newcommand{\go}{\stackrel {\rm G}{\leq}}     
\newcommand{\geg}{\stackrel {\rm G}{\geq}}
\newcommand{\gegs}{\stackrel {\rm G}{>}}
\newcommand{\gos}{\stackrel {\rm G}{<}}
\newcommand{\diar}{{\sr \swarrow}}
\newcommand{\rar}{\rightarrow}
\newtheorem*{theorem}{Theorem}
\newtheorem*{lemma}{Lemma}

\newtheorem*{prop}{Proposition}

\begin{document}
\title[Descendants of Richardson component]
{\bf On orbital variety closures in ${\gs\gl}_n$\\
\bf II.  Descendants of a Richardson orbital variety} 
\author{Anna Melnikov}\thanks{This work
was partially supported by the EEC program TMR-grant ERB FMRX-C
T97-0100}
\address{Department of Mathematics,
University of Haifa,
Haifa 31905, Israel}
\email{melnikov@math.haifa.ac.il}
\begin{abstract}
For  a semisimple Lie algebra $\gog$ the orbit
method attempts to assign representations
of $\gog$ to (coadjoint) orbits in $\gog^*.$  Orbital
varieties are particular Lagrangian subvarieties of such orbits leading 
to highest weight representations of $\gog.$ 
In ${\gs\gl}_n$ orbital varieties are described by 
Young tableaux. 
In this paper we consider so called Richardson orbital varieties in $\gs\gl_n.$
A Richardson orbital variety is an orbital variety
whose closure is a standard nilradical. We show that in $\gs\gl_n$ a Richardson
orbital variety closure is a union of orbital varieties. We give a complete
combinatorial description of such closures in terms of Young tableaux.
\par
This is the second paper in the series of three papers devoted to a combinatorial
description of orbital variety closures in $\gs\gl_n$ in terms of Young tableaux.
\end{abstract}
\maketitle

\section {\bf Introduction}
\subsection{}\label{1.1}
This is a continuation of Part I 
(\cite{M}) whose notation we retain. Our main goal in this series of papers is
to give a description of an orbital variety closure in ${\gs\gl}_n.$
In this paper we give a complete description of the closure of a so called Richardson
orbital variety. 
\par
For the convenience of the reader we repeat necessary notation and results from Part I
which can be formulated in short. If the formulation is too long, as for example, in the
case of Robinson-Schensted procedure we provide the exact reference to the subsection 
of Part I.
\subsection{}\label{1.2} 
Let $\bG$ be a connected semisimple 
finite dimensional
complex algebraic group. Let $\gog$ be its Lie algebra.
Fix some triangular decomposition  $\gog=\gn\bigoplus
\gh\bigoplus\gn^-.$  Take some $x\in\gn$ and consider $\Oscr=\Oscr_x$ a nilpotent orbit of $\gog$
under the adjoint action of $\bG$ (that is $\Oscr=\bG x=\{gxg^{-1}\ |\ g\in\bG\}$).
Since $\gog$ is semi-simple we can identify $\gog^*$ with $\gog$ through the Killing form.
This identification gives an adjoint orbit a simplectic structure. 
By the results of R. Steinberg, N. Spaltenstein and A. Joseph  (cf. ~\cite{J}) $\Oscr\cap\gn$ 
is equidimentional and Lagrangian.
Its irreducible components are called orbital varieties
associated to $\Oscr.$
\par
Given an orbital variety $\Vscr$ we denote by $\Oscr_\Vscr$ the 
orbit $\Vscr$ is associated to.
\subsection{}\label{1.3}
Let $R\st\gh^*$ denote the set of non-zero roots,  
$R^+$ the set of positive roots corresponding to $\gn$ and 
$\Pi\st R^+$ the resulting set of simple roots. Let $X_{\al}:={\Co} x_{\al}$ denote the root subspace 
corresponding to $\al\in R.$ Then $\gn=\bigoplus\limits_{\al\in R^+}X_{\al}.$ Set 
$\gn\cap^w\gn:=\bigoplus\limits_{\al\in R^+\cap^w R^+}X_{\al}.$
Let $\bB$ be the standard
Borel subgroup of $\bf G$, i.e.  such that $\rm {Lie}({\bf B})=\gb=\gh\bigoplus\gn.$
$\bf B$ acts by conjugation on $\mathfrak n$ and its subsets.
Let $W$  be the Weyl group for 
the pair $(\gog,\ \gh).$ By Steinberg's construction there
exists a surjection from $W$ onto the set of orbital varieties
defined by $w\rar \ov{\bB(\gn\cap^w \gn)}=:\ov{\Vscr}_w$ (cf. Part I, 2.1.2, 2.1.3).
\par
Note that there exists the unique nilpotent orbit $\Oscr_w$
such that $\ov\Oscr_w=\ov{\bG(\gn\cap^w \gn)}.$ Obviously one has
$\Oscr_w=\Oscr_{\Vscr_w}.$ 
\subsection{}\label{1.4}
Take $\Iscr\st \Pi,$ let ${\bP}_{\Iscr}$ denote the unique 
standard parabolic subgroup of $\bG$ generated by
the standard parabolic subgroups
${\bP}_{\al}\ :\ \al\in \Iscr.$ Let $\bM_{\Iscr}$ be the unipotent 
radical of $\bP_{\Iscr}$ and $\gm_{\sr \Iscr}=\Lie({\bM}_{\Iscr})$
be the corresponding nilradical in $\gn.$
\par
Let $W_{\Iscr}:=
<s_{\al}\ :\ \al\in \Iscr>$ be the corresponding parabolic subgroup of $W.$
Let $w_{\sr \Iscr}$ be the unique longest element of $W_{\Iscr},$
that is such that
$$w_{\sr \Iscr}(\al)\in\begin{cases} R^- & {\rm if}\ \al\in\Iscr\cr
                             R^+ & {\rm if}\ \al\in R^+\backslash \sum\limits_{\be\in\Iscr}{\Na}
\cdot\be \cr\end{cases}\eqno{(*)}$$ 
Note that $\gn\cap^{w_{\sr \Iscr}}\gn=\gm_{\sr \Iscr}.$ Thus
$\ov{\Vscr}_{w_\Iscr}=\ov{\bB(\gn\cap^{w_\Iscr} \gn)}=
\gm_{\sr \Iscr}.$ Just to simplify the notation we denote 
$\Vscr_\Iscr:=\Vscr_{w_\Iscr}.$ 
\par
Set $\Oscr_\Iscr:=\Oscr_{\Vscr_\Iscr}.$ One calls
$\Vscr_\Iscr$ a Richardson orbital variety
(associated to $\Oscr_\Iscr$) defined by $\Iscr.$
\subsection{}\label{1.5}
From now on we consider only
the case of $\gog=\gs\gl_n.$ Every nilpotent orbit in $\gog$
is defined by its Jordan form which in turn is completely defined
by the set of lengths of its Jordan blocks. So we get a bijection
$\varphi$
from the set of partitions of $n$ 
onto the set of nilpotent orbits in $\gog.$ 
We write a partition 
in decreasing order that is $\lam=(\lam_1,\lam_2,\cdots,\lam_k)$ where $\sum_{i=1}^k\lam_i=n$
and $\lam_1\geq\lam_2\geq\cdots\geq l_k> 0.$
We set $\Oscr_\lam:=\varphi(\lam).$ 
\par
Define a partial order on partitions of $n$ as follows. 
Let $\lam=(\lam_{\sr 1},\cdots,\lam_k)$ and 
$\mu = (\mu_{\sr 1},\cdots,\mu_j) $
be partitions of $n.$ If $j\ne k$ complete the partition
with less number of elements by $0$ so that we can consider 
that both partitions have $\max(j,k)$ elements.
Define $\lam\leq \mu$ if for each
$i\ :\ 1\leq i\leq \max(j,k)$ one has
$$ \sum\limits^i_{m=1}\lam_m \geq \sum\limits^i_{m=1}\mu_m. $$
\par
Then by the result of Gerstenhaber (cf. Part I 2.3.2)
$$\ov\Oscr_\lam=\coprod\limits_{\mu\geq \lam}\Oscr_\mu.$$
We define a partial order on nilpotent orbits by 
$\Oscr\leq\Oscr\pr$ if $\Oscr\pr\st\ov\Oscr.$ We say that 
$\Oscr\pr\ :\ \ov\Oscr\pr\subsetneq\ov\Oscr$
is a descendant of $\Oscr$ if for any $\Oscr\prpr$ such that 
$\Oscr\leq\Oscr\prpr\leq \Oscr\pr$ one has $\Oscr\prpr=\Oscr$
or $\Oscr\prpr=\Oscr\pr.$
\subsection{}\label{1.6} 
Recall the notation and the
conventions of Part I, 2.2.1. In particular recall that 
in our case $\Pi=\{\al_i\}_{i=1}^{n-1}.$   
Recall from Part I, 2.4.5, 2.4.6 the combinatorial 
characterization of orbital varieties in 
$\gog$ in terms of Young tableaux. It is constructed as follows.
As we explained in \ref{1.3} by Steinberg construction one has
a surjection from Weyl group $W$ onto the set of orbital varieties.
Set $\Cscr_{\Vscr}:=\{w\in W\ |\ \ov\Vscr=\ov{\Vscr}_w\}.$ We get
a partition of Weyl group into so called geometric cells.
We also define $\Cscr_w:=\Cscr_{\Vscr_w}.$
\par
Identify Weyl group of $\gs\gl_n$ with symmetric
group $\bS_n$ (cf. Part I, 2.2.2).
Robinson-Schensted procedure (cf. Part I, 2.4.6) gives a
surjection $w\mapsto T(w)$ from $\bS_n$ onto $\bT_n$ the set  of standard 
Young tableaux with $n$ entries.  By R. Steinberg (cf. Part I, 2.4.6)
one has $y\in \Cscr_w$ if and only if $T(y)=T(w).$ In such a way we
we get a bijection $\phi(T(w)):=\Vscr_w$ from the set
of standard Young tableaux onto the set of orbital varieties.
Given a standard Young tableau $T$ set ${\Vscr}_T:=\phi(T).$
Respectively given an orbital variety $\Vscr$ set $T_{\Vscr}:=\phi^{\sr -1}(\Vscr).$
\par
Given a Young tableau $T$ let $\sh(T)=\lam$ be the  partition of $n$ from which $T$
was built, that is $\lambda=(\lambda_1,\ldots,\lambda_k)$ where $\lambda_i$ is the length
or $i$-th row of $T.$ Then $\Vscr_T$ is associated to 
$\Oscr_{\lam}$ so that
the characterization of orbital varieties by Young tableaux is compatible
with the characterization of nilpotent orbits by partitions.
\subsection{}\label{1.6a}
We give a combinatorial description of a Richardson orbital variety closure
in terms of Young tableaux in the spirit of Gerstenhaber's
construction.
\par
Given a standard Young tableau $T\in \bT_n$ and a natural number $u\ :\ 1\leq u\leq n$
let  
$r_{\sr T}(u)$ denote the number of the row of $T$ $u$ belongs to.
Recall from Part I, 2.4.14 that $\tau(T)=\{\al_i\ :\ r_{\sr T}(i+1)>r_{\sr T}(i)\}.$ 
Recall also from there that $\tau({\Vscr}_T)=\tau(T).$ 

We show in \ref{2.1} that
$$\ov\Vscr_\Iscr=\bigcup_{\tau(\Vscr)\supseteq\Iscr}\Vscr.\eqno{(**)}$$
This formula resembles Gerstenhaber's
result in its simplicity and clearness.
\subsection{}\label{1.7} We define a 
partial order on the set of orbital varieties by $\Vscr\go\Wscr$ if $\Wscr\st\ov\Vscr.$
We call it the geometric order.
This order is compatible with the order on nilpotent orbits defined in \ref{1.5} in the sense
that $\Vscr\go\Wscr$ implies that $\Oscr_\Vscr\leq \Oscr_\Wscr.$ Moreover by Part I, 4.1.8
for any orbital variety $\Vscr$ in $\gs\gl_n$ and any
nilpotent orbit $\Oscr$ such that $\Oscr>\Oscr_\Vscr$ there exist an orbital variety $\Wscr$ associated to
$\Oscr$ such that $\Vscr\gos \Wscr.$
\par
We say that $\Wscr\ :\ \Wscr\gegs\Vscr$ is a geometric descendant 
of $\Vscr$
if for any $\Yscr$ such that 
$\Vscr\go\Yscr\go\Wscr$ one has $\Yscr=\Vscr$ or $\Yscr=\Wscr.$
\par
The description of the set of descendants of a nilpotent orbit
is obtained as an easy corollary of Gerstenhaber's
formula (cf. Part I, 2.3.3). However
the description of the set of geometric descendants
of a Richardson orbital variety  cannot be easily obtained from
\ref{1.6a} $(**).$
Yet the result is very clear-cut. This is the main result of the paper. 
Its formulation
demands some additional combinatorial notation so we refer the reader
to \ref{2.6} for the exact statement of the result.
\subsection{}\label{1.8} By ~\cite{Sp} $\gn\cap\Oscr$ is a union of orbital varieties
associated to $\Oscr.$ Moreover 
by ~\cite{Me} 
$$\ov{\gn\cap\Oscr}=\gn\cap\ov\Oscr.$$
As we show in \ref{2.1} for any $\Oscr\st\ov{\Oscr}_{\Iscr}$
one has $\gm_{\Iscr}\cap\Oscr$ is a union of all orbital varieties
(associated to $\Oscr$) whose $\tau$-invariant
contains $\Iscr.$ A natural question is whether for such $\Oscr$
one has as well that
$$\ov{\gm_{\sr \Iscr}\cap\Oscr}=\gm_{\sr \Iscr}\cap\ov\Oscr. \eqno{(***)}$$ 
\par
Theorem \ref{2.6} provides
in particular that $\Vscr$ being a descendant of $\Vscr_\Iscr$
does not necessarily imply that $\Oscr_\Vscr$ is a descendant 
of $\Oscr_\Iscr.$ Using this fact we show in \ref{4.1} that in general 
equality $(***)$ does not hold.
\subsection{}\label{1.9} 
In ~\cite{JM} the ideals of definition of an
orbital variety closure of codimension 1 in some nilradical were 
constructed as well as strong quantization of such orbital 
variety. The description of the set of geometric descendants
of a Richardson orbital variety is central for the generalization of 
the results of ~\cite{JM} to orbital varieties of codimension greater than 
1 in a nilradical.
\subsection{}\label{1.10}
A natural and interesting problem 
is to determine which orbital variety closures are complete intersections. This 
is very important in particular for strong quantization of an orbital variety. 
It is more or less obvious that
the majority of orbital varieties are not complete intersections.
Concentrating on the descendants of a Richardson orbital variety
one can see at once that all the descendants of codimension 1 in $\gm_{\sr \Iscr}$
are complete intersections just by Krull theorem. What can be said about descendants
of a Richardson orbital variety of  codimension greater than 1 in $\gm_{\sr \Iscr}$?
It is very easy to check that in $\gs\gl_n$ 
for $n\leq 5$ all the descendants of a Richardson orbital variety are complete
intersections. 
We give an example in $\gs\gl_{\sr 6}$  of $\Vscr_\Iscr$ and its descendant $\Wscr$ 
such that 
\begin{itemize}
\item[(i)] $\codim_{\gm_\Iscr}\ov\Wscr=2;$
\item[(ii)] $\Oscr_\Wscr$ is a descendant of $\Oscr_\Iscr;$
\item[(iii)] $\Wscr$ is not a complete intersection.
\end{itemize}
We construct the ideal
of definition for this $\Wscr.$ 
\par
This example shows that even among the descendants of a Richardson orbital variety
the majority of those of codimension greater than 1 in $\gm_{\sr \Iscr}$ are not
complete intersections. 
\subsection{}\label{1.11}
The body of the paper consists of three sections.
\par
In Section 2 we develop the combinatoric notation essential to state 
the main theorem (formulated in \ref{2.6}).
\par
Section 3 is devoted to the proof of theorem \ref{2.6}. Here we develop 
the combinatorics connected to the non-trivial involution
$\psi$ of  the Dynkin diagram of ${\gs\gl}_n.$ This involution induces
involution on Weyl group $W.$ 
Given $w\in W$ in word form
we give a formula for $\psi(w)$ in a word form in \ref{3.7}. 
As well involution $\psi$ induces order preserving involution on the set
of Young tableaux. In general there is no straightforward formula determining 
$\psi(T).$ However we develop such formulas for the tableaux considered in the paper
and use them as one of the central tools in our proof.
\par
Finally in section 4 we discuss in detail \ref{1.8} and \ref{1.10}.
\par
In the end one can find the index of notation in which 
symbols appearing frequently are given with the subsection
where they are defined. 
%
\section{\bf A Richardson orbital variety closure and its descendants}
\subsection{}\label{2.0} 
Let us recall the notation and the results concerning
Weyl group and in particular $\bS_n$ from Part I that we use
in this section.
\par
For $w\in W$
set $S(w):=R^+\cap^w R^-.$ 
The map $w\mapsto S(w)$ is injective, so we define a partial order relation
on $W$ by taking $y\dor w$ if $S(y)\subseteq S(w).$ 
It is called the Duflo order.
\par
Recall notation $X_{\al}$ from \ref{1.3}.
Since $\gn\cap^w\gn=\bigoplus\limits_{\al\in R^+\cap^w R^+}X_{\al}$
one has $y\dor w$ if and only if $\gn\cap^w\gn\subseteq \gn\cap^y\gn.$
In particular if $y\dor w$ then $\Vscr_w\go\Vscr_y.$
\par
Set $\tau(w):=S(w)\cap\Pi.$ Recall notion  of the standard Borel subalgebra $\gb$
from \ref{1.3}.
For $\al\in \Pi$ let $\bP_{\al}$ be the standard parabolic
subgroup such that $\Lie(\bP_{\al})=\gb\oplus X_{-\al}.$ Given an orbital variety
$\Vscr$ set $\bP_{\al}(\Vscr):=\{gXg^{\sr -1}\ |\ X\in \Vscr,\ g\in \bP_{\al}\}$ 
and let $\tau(\Vscr):=\{\al\in \Pi\ |\ \bP_{\al}(\Vscr)=\Vscr\}.$
\par 
Let us consider $\gog=\gs\gl_n.$ In that case $\Pi=\{\al_i\}_{i=1}^{n-1}$
and any $\al\in R^+$ is $\al_{i,j}=\sum\limits_{k=i}^{j-1}\al_k$ for some
$i,j\ :\ 1\leq i<j\leq n.$ As well in that case 
$W$ is identified with $\bS_n$ (cf. Part I, 2.2.1, 2.2.2). 
We write $w\in \bS_n$ in a word form $w=[a_1,\ldots,a_n]$ where $a_i=w(i).$
Put $p_w(i)=j$ if $i=a_j$, that is $p_w(i)$ is its place (index) in word
$w.$ By  Part I, 2.2.4 one has $\al_{i,j}\in S(w)$
iff $p_w(i)>p_w(j).$
\par
Recall $\tau(T)$ from \ref{1.6a}. By Part I, 2.1.7, 2.4.14 $\tau(w)=\tau(\Vscr_w)=\tau(T(w)).$
\par
We also need the following notation from Part I, 2.2.5. Given a word
$w=[a_1,\ldots,a_n]$ let $<w>:=\{a_i\}_{i=1}^n$ be the set of its entries. 
\begin{itemize}
\item[(i)] Set $\ov w:=[a_n,\ldots,a_1]$ to be the word with order
reverse to the order of $w$.
\item[(ii)] Given words $x=[a_1,\ldots, a_k],\ y=[b_1,\ldots, b_m]$ such that
$<x>\cap <y>=\emptyset$ set 
$[x,y]=[a_1,\ldots,a_k,b_1,\ldots,b_m]$ to be
a colligation.
\end{itemize}
Given a fixed set $E$
of $n$ distinct positive integers we let ${\boS}_n$ or ${\boS}_E$ denote the
set of words $w$ such that the set of its entries $<w>=E.$ 
\subsection{}\label{2.01}
As well we need to recall some notation concerning Young tableaux.
Given a Young tableau $T$ recall from Part I, 2.4.2 that $T^i_j$ is the entry on the intersection
of $i$-th row and $j$-th column and $\omega^i(T)$ is the last`< entry of row $i$ of $T.$ 
As well $T_j$ is $j$-th column of $T;$\ \ $T^i$
is $i$-th row of $T;$ notation $T^{i,j}$ means subtableau
of $T$ containing all the rows from $i$-th row to $j$-th row  and finally $T^{i,\infty}$ 
is the subtableau of $T$ containing all the rows from $i$-th row and down.
Put $|T_j|$ (resp. $|T^i|$) to be the length of $j$-th column (resp. $i$-th row)
of $T.$ 
\par
For example if 
$$T=
\vcenter{
\halign{& \hfill#\hfill
\tabskip4pt\cr
\multispan{11}{\hrulefill}\cr
\ssa
\vb & 1 &  &  2 & &  4 && 8 && 9 &\ts\vb\cr
\vsa
&&&&&&\multispan{5}{\hrulefill}\cr
\ssa
\vb & 3 &  &  5 && 6 & \ts\vb\cr
\vsa
&&\multispan{5}{\hrulefill}\cr
\ssa
\vb & 7 & \ts\vb\cr
\vsa
&&\cr
\ssa
\vb & 10 & \ts\vb\cr
\vsa
\multispan{3}{\hrulefill}\cr}}$$
then $T^2_2=5,\ T^1_4=8;\ \omega^1(T)=9,\ \omega^2(T)=6,\ \omega^3(T)=7,\ \omega^4(T)=10$ and
$$T_3=
\vcenter{
\halign{& \hfill#\hfill
\tabskip4pt\cr
\multispan{3}{\hrulefill}\cr
\ssa
\vb & 4 &\ts\vb\cr
\vsa
&&\cr
\ssa
\vb & 6 & \ts\vb\cr
\vsa
\multispan{3}{\hrulefill}\cr}},\ 
 T^2=
\vcenter{
\halign{& \hfill#\hfill
\tabskip4pt\cr
\multispan{7}{\hrulefill}\cr
\ssa
\vb & 3 &  &  5 && 6 & \ts\vb\cr
\vsa
\multispan{7}{\hrulefill}\cr}},\ T^{2,3}=
\vcenter{
\halign{& \hfill#\hfill
\tabskip4pt\cr
\multispan{7}{\hrulefill}\cr
\ssa
\vb & 3 &  &  5 && 6 & \ts\vb\cr
\vsa
&&\multispan{5}{\hrulefill}\cr
\ssa
\vb & 7 & \ts\vb\cr
\vsa
\multispan{3}{\hrulefill}\cr}},\ T^{2,\infty}=
\vcenter{
\halign{& \hfill#\hfill
\tabskip4pt\cr
\multispan{7}{\hrulefill}\cr
\ssa
\vb & 3 &  &  5 && 6 & \ts\vb\cr
\vsa
&&\multispan{5}{\hrulefill}\cr
\ssa
\vb & 7 & \ts\vb\cr
\vsa
&&\cr
\ssa
\vb & 10 & \ts\vb\cr
\vsa
\multispan{3}{\hrulefill}\cr}}$$
\par
Given a tableau $T$ let 
$<T>$ denote the set of its entries.
Given a fixed set $E$ of $n$ distinct positive integers, let ${\boT}_n$ or ${\boT}_E$
be the set of Young tableaux $T$ such that $<T>=E.$
\par
Recall from Part I, 2.4.3 the definition of a tableau $(T,S)$ where $<T>\cap<S>=\emptyset.$ Note that if
for every row $i: 1\leq i\leq \min(|S_1|,|T_1|)$ one has $\omega^i(T)<S^i_1$ then $(T,S)$ is obtained
by shifting  cells of $S$ to the left.
\par 
Recall also from Part I, 2.4.2 that given a tableau $T$ we set $T^{\dagger}$ to be a transposed tableau.
For example if $T$ is the tableau above then
$$T^{\dagger}=
 \vcenter{
\halign{& \hfill#\hfill
\tabskip4pt\cr
\multispan{9}{\hrulefill}\cr
\ssa
\vb & 1 &  &  3 & &  7 && 10 &\ts\vb\cr
\vsa
&&&&&&\multispan{3}{\hrulefill}\cr
\ssa
\vb & 2 &  &  5 && 6 & \ts\vb\cr
\vsa
&&&&\multispan{3}{\hrulefill}\cr
\ssa
\vb & 4 & & 6 &\ts\vb\cr
\vsa
&&\multispan{3}{\hrulefill}\cr
\ssa
\vb & 8& \ts\vb\cr
\vsa
&&\cr
\ssa
\vb & 9& \ts\vb\cr
\vsa
\multispan{3}{\hrulefill}\cr}}$$

\subsection{}\label{2.1}
We begin with the result $(**)$ from \ref{1.6a}. This is a
straightforward corollary of ~\cite[2.13]{JM}. Here we prove it for completeness.
\begin{theorem} For any $\Iscr\st\Pi$ one has
$$\ov{\Vscr}_{\Iscr}=\bigcup_{\tau(\Vscr)\supseteq\Iscr}\Vscr.$$
\end{theorem}
\Pf
First of all note that  for any $w\in W$ one has
$w\dg w_{\sr \Iscr}$ if and only if 
$\tau(w)\supseteq \Iscr.$ Thus
$$\ov{\Vscr}_{\Iscr} \supseteq \bigcup\limits_{\tau(\Vscr)\supseteq \Iscr}{\Vscr}.$$
Moreover since $\ov{{\bB}(\gn\cap^{w_\Iscr}\gn)}=\gn\cap^{w_\Iscr}\gn$
one has $\Vscr\geg\Vscr_\Iscr$ if and only if $\Vscr\dg\Vscr_\Iscr.$
Hence $\Vscr\geg\Vscr_\Iscr$ if and only if $\tau(\Vscr)\supseteq \Iscr.$
\par
On the other hand
$$\ov{\Vscr}_{\Iscr}={\gm}_{\sr \Iscr}=
\coprod_{\Oscr\geq{\Oscr}_{\Iscr}}\Oscr\cap\gm_{\Iscr}$$
So to complete the proof we must show that for any 
$\Oscr>{\Oscr}_{\Iscr}$
any irreducible component
of $\gm_{\Iscr}\cap\Oscr$ is an orbital variety.
By  ~\cite[p. 456, last corollary]{Sp1} $\gm_{\Iscr}\cap\Oscr$
is equidimensional. 
By Part I, 4.1.8 it contains at least one orbital
variety and so for every irreducible 
component $\Vscr$ of $\Oscr\cap\gm_{\Iscr}$ one has
$\dim \Vscr=\frac{1}{2}\, \dim \Oscr.$  
Yet $\Vscr\st \Oscr\cap\gn,$
which also has dimension $\frac{1}{2}\, \dim \Oscr.$ 
Hence $\Vscr$ is an orbital
variety associated to $\Oscr.$ 
\QED
\subsection{}\label{2.2} Given a tableau $T$ recall from \ref{1.6a} that
$r_{\sr T}(u)$ denotes the number of the row $u$ belongs to.
\par
Given $\Iscr\st\Pi.$ Set $T_\Iscr:=T_{\Vscr_\Iscr}.$
Since $T_\Iscr$ is the minimal (both in the geometric and the Duflo
orders) tableau such that $\tau(T)=\Iscr$ one has
$$r_{\sr T_{\Iscr}}(i)=\begin{cases} 
1, & {\rm if}\ i=1\ {\rm or}\ \al_{i-1}\not\in\Iscr,\cr
         r_{\sr T_\Iscr}(i-1)+1, & {\rm if}\ \al_{i-1}\in\Iscr.\cr\end{cases}$$
A useful way to present $T_\Iscr$ is provided in ~\cite[2.12]{JM}.
Partition
$\{1,2,\ldots,n\}$ into connected subsets 
$C_j:=\{b_j,b_j+1,\ldots,b_{j+1}-1\}$ by choosing a strictly 
increasing sequence $1=b_1<b_2<\cdots<b_{l+1}=n+1.$
Setting $\Iscr=\{\al_i\ |\ i,i+1\ {\rm belong\ to\ some\ }C_j\}$
defines a bijection between the set of all such  partitions and 
the set of subsets of $\Pi.$ Given $\Iscr\st \Pi,$ let
$\{C_i^{\Iscr}\ :\ i=1,2,\ldots,l\}$ be the corresponding 
connected subsets which we view as columns. (Sometimes we may
omit the $\Iscr$ superscript.) Then
in the notation of \ref{2.01} we 
have $T_{\Iscr}=(C^{\Iscr}_1,C^{\Iscr}_2,\ldots,C^{\Iscr}_l).$
Of course this involves some sliding of boxes to the left.
However there are some advantages in this presentation. For example
$T_\Iscr -1$ is obtained by simply replacing $C^{\Iscr}_1$ by $C^{\Iscr}_1-1$ and
$T_\Iscr-n$ is obtained by simply replacing $C^{\Iscr}_l$ by $C^{\Iscr}_l-n.$
\par
We call $C^\Iscr_i$ the $i-$th chain of $T_\Iscr$ and 
$T_\Iscr=(C^{\Iscr}_1, C^{\Iscr}_2,\ldots, C^{\Iscr}_l)$ the chain form of
$T_\Iscr.$
\par
Set $c_i:=|C_i|$ and $\si_i:=C^{c_i}_i$ which is the
largest entry in $C_i.$ 
\par
For example
consider $\Iscr=\{\al_1,\al_4,\al_5,\al_6,\al_9\}$ in $\gs\gl_{10}.$ Then
$$T_\Iscr =
\vcenter{
\halign{& \hfill#\hfill
\tabskip4pt\cr
\multispan{11}{\hrulefill}\cr
\ssa
\vb & 1 &  &  3 & &  4 && 8 && 9 &\ts\vb\cr
\vsa
&&&&&&\multispan{5}{\hrulefill}\cr
\ssa
\vb & 2 &  &  5 && 10 & \ts\vb\cr
\vsa
&&\multispan{5}{\hrulefill}\cr
\ssa
\vb & 6 & \ts\vb\cr
\vsa
&&\cr
\ssa
\vb & 7 & \ts\vb\cr
\vsa
\multispan{3}{\hrulefill}\cr}}$$
or in chain form $T_\Iscr=(C_1,C_2,C_3,C_4,C_5)$ where $C_1=\{1,2\},$\ $C_2=\{3\},$\ 
$C_3=\{4,5,6,7\},$\ $C_4=\{8\}$ and $C_5=\{9,10\}.$ Here one has $c_{\sr 1}=c_{\sr 5}=2,\ 
c_{\sr 2}=c_{\sr 4}=1$ and $c_{\sr 3}=4.$ As well $\si_{\sr 1}=2,\ \si_{\sr 2}=3,\ 
 \si_{\sr 3}=7,\ \si_{\sr 4}=8,\ \si_{\sr 5}=10.$
\subsection{}\label{2.3}
Recall the notion $w_{\sr \Iscr}$ from \ref{1.4}. 
By \ref{2.0} one has $\al_i\in \tau(w_{\sr \Iscr})$ if and only if 
$p_{w_{\Iscr}}(i)>p_{w_{\Iscr}}(i+1)$  and by \ref{1.4} $(*)$ we get
$\al_i\in \tau(w_{\sr \Iscr})$ if and only if $\al_i\in\Iscr.$ Thus
$p_{w_{\Iscr}}(i)>p_{w_{\Iscr}}(i+1)$ if and only if $\al_i\in\Iscr$
which can be formulated as
\begin{prop}[ {~\cite[2.12]{JM}} ] For all $\Iscr\st \Pi,$ the word form
of $w_{\sr \Iscr}$ is given by 
$$w_{\sr \Iscr}=[\ov{C^{\Iscr}_1},\ov{C^{\Iscr}_2},\ldots,
\ov{C^{\Iscr}_l}].$$
\end{prop}
\subsection{}\label{2.4} For $s,t\ : 1\leq s<t\leq l$ set 
$\hat T_{s,t}:=(C_s,C_{s+1},\ldots,C_t).$ 

For $i:\ 2\leq i\leq l$ and for some $j:\ 1\leq j\leq c_i$ set $C_i^{1,j-1}:=(C_i^1,C_i^2,\ldots, C_i^{j-1})$ and 
$C_i(j):=(C^j_i,\ldots,C^{c_i}_i).$ We call $C_i(j)$ a $j-$tail of the chain $C_i.$
If there exist $s<i$ such that $c_s\geq j$ then set 
$$T\pr=\left(\begin{array}{c}(\hat T_{1,i-1})^{1,c_s},C_i^{1,j-1}\cr
             (\hat T_{1,i-1})^{c_s+1,\infty},C_i(j)\cr\end{array}\right)\quad {\rm and}\quad T_\Iscr(c_s\diar C_i(j))=(T\pr,\hat T_{i+1,l}). \eqno{(*)}$$
Note that $T\pr$ is obtained from $\hat T_{1,i}$ just by moving 
$j-$tail of the last chain down to the
part $(\hat T_{1,i-1})^{c_s+1,\infty}.$ Since $|(\hat T_{1,i-1})^{c_s+1}|<|(\hat T_{1,i-1})^{c_s}|$
we get that $T\pr$ and 
respectively $T_\Iscr(c_s\diar C_i(j))$ are Young tableaux.
\par
For example let $\Iscr=\{\al_1,\al_3,\al_4,\al_5,\al_7,\al_8\}\st \Pi_{1,10}$ that is
$$ T_\Iscr =
\vcenter{
\halign{& \hfill#\hfill
\tabskip4pt\cr
\multispan{9}{\hrulefill}\cr
\ssa
\vb & 1 &  &  3 & &  7 && 10& \ts\vb\cr
\vsa
&&&&&&\multispan{3}{\hrulefill}\cr
\ssa
\vb & 2 &  &  4 & & 8 & \ts\vb\cr
\vsa
&&&&\multispan{3}{\hrulefill}\cr
\ssa
\vb & 5 && 9 & \ts\vb\cr
\vsa
&&\multispan{3}{\hrulefill}\cr
\ssa
\vb & 6 & \ts\vb\cr
\vsa
\multispan{3}{\hrulefill}\cr
}}$$
Then 
$$T_\Iscr(2\diar C_2(2)) =
\vcenter{
\halign{& \hfill#\hfill
\tabskip4pt\cr
\multispan{9}{\hrulefill}\cr
\ssa
\vb & 1 &  &  3 & &  7 && 10& \ts\vb\cr
\vsa
&&&&\multispan{5}{\hrulefill}\cr
\ssa
\vb & 2 &  & 8 & \ts\vb\cr
\vsa
\multispan{3}{\hrulefill}&&&&\cr
\ssa
\vb & 4 &\vb & 9 &\ts\vb\cr
\vsa
&&\multispan{3}{\hrulefill}\cr
\ssa
\vb & 5 & \ts\vb\cr
\vsa
&&\cr
\ssa
\vb & 6 & \ts\vb\cr
\vsa
\multispan{3}{\hrulefill}\cr
}};\quad
T_\Iscr(2\diar C_3(2)) =
\vcenter{
\halign{& \hfill#\hfill
\tabskip4pt\cr
\multispan{9}{\hrulefill}\cr
\ssa
\vb & 1 &  &  3 & &  7 && 10& \ts\vb\cr
\vsa
&& && \multispan{5}{\hrulefill}\cr
\ssa
\vb & 2 &  &  4 & \ts\vb\cr
\vsa
&&\multispan{3}{\hrulefill}\cr
\ssa
\vb & 5 &\vb& 8& \ts\vb\cr
\vsa
&&&&\cr
\ssa
\vb & 6 &\vb& 9 &\ts\vb\cr
\vsa
\multispan{5}{\hrulefill}\cr
}}; 
$$
$$T_\Iscr(3\diar C_4(1))=
\vcenter{
\halign{& \hfill#\hfill
\tabskip4pt\cr
\multispan{7}{\hrulefill}\cr
\ssa
\vb & 1 &  &  3 & &  7 & \ts\vb\cr
\vsa
&&&&&&\multispan{0}{\hrulefill}\cr
\ssa
\vb & 2 &  &  4 & & 8 & \ts\vb\cr
\vsa
&&&&\multispan{3}{\hrulefill}\cr
\ssa
\vb & 5 && 9 & \ts\vb\cr
\vsa
&&\multispan{3}{\hrulefill}\cr
\ssa
\vb & 6 &\vb& 10& \ts\vb\cr
\vsa
\multispan{5}{\hrulefill}\cr
}}.
$$

\subsection{}\label{2.5}
Consider $\hat T_{1,i-1},$ that is, look to the left of chain $C_i.$ Set
$$\begin{array}{rcl}\next_l(i)&:=&\begin{cases}\min\{c_j\ |\ j<i,\ c_j\geq c_i\},& 
{\rm if\ there\ exists}\ j<i\ {\rm such\ that}\ c_j\geq c_i;\cr
                     0, & {\rm otherwise}.\cr\end{cases}\\  
\prev_l(i)&:=&\begin{cases}\max\{c_j\ |\ j<i,\ c_j\leq c_i\},& {\rm if\ there\ exists}\ j<i\ 
{\rm such\ that}\ c_j\leq c_i;\cr
                     n,& {\rm otherwise}.\cr\end{cases}\\
\sprev_l(i)&:=&\begin{cases}\max\{c_j\ |\ j<i,\ c_j<c_i\},& {\rm if\ there\ exists}\ j<i\ 
{\rm such\ that}\ c_j<c_i;\cr
                     n,& {\rm otherwise}.\cr\end{cases}\\
\end{array}.$$
Consider $\hat T_{i+1,l}$ that is look to the right of chain $C_i.$ Set
$$\begin{array}{rcl}
\snext_r(i)&:=&\begin{cases}\min\{c_j\ |\ j>i,\ c_j>c_i\},& {\rm if\ there\ exists}\ j>i
{\rm such\ that}\ c_j> c_i;\cr
                     n, & {\rm otherwise}.\cr\end{cases}\\
\sprev_r(i)&:=&\begin{cases}\max\{c_j\ |\ j>i,\ c_j< c_i\},& {\rm if\ there\ exists}\ j>i 
{\rm such\ that}\ c_j<c_i;\cr
                    0,& {\rm otherwise.}\cr\end{cases}\\
\end{array}$$
We use the nomenclature ``sprev'' (resp. ``snext'') to emphasize that this is a strictly
previous (resp. a strictly next) number.
\par
Set
$$T_\Iscr(i,\next):=\begin{cases}T_\Iscr(\next_l(i)\diar C_i(c_i)), & {\rm if}\ \next_l(i)\ne 0,\cr 
                          \emptyset, & {\rm otherwise}.\cr\end{cases}$$
That is consider $\hat T_{1,i}$ and move the box with $\si_i,$ 
which is the largest number of  $\hat T_{1,i},$ 
down to the first possible row (if such row exists), then add 
to the new tableau $\hat T_{i+1,l}.$
\parno
Set
$$T_\Iscr(i,\prev)=\begin{cases}T_\Iscr(\sprev_l(i)\diar C_i(\sprev_l(i))), & {\rm if}\ 
\prev_l(i)=\sprev_l(i)\ne n,\cr
                          \emptyset, & {\rm otherwise}.\cr\end{cases}$$
That is consider $\hat T_{1,i}$ and move the smallest possible
tail of the last column one row down (if this is possible), then add 
to the new tableau $\hat T_{i+1,l}.$
\parno
Set $T_\Iscr(i):=\{T_\Iscr(i,\next),T_\Iscr(i,\prev)\}.$
\par
For example consider $T_\Iscr$ from \ref{2.4}. Then
$$ T_\Iscr(2)=\left\{
\vcenter{
\halign{& \hfill#\hfill
\tabskip4pt\cr
\multispan{9}{\hrulefill}\cr
\ssa
\vb & 1 &  &  3 & &  7 & & 10 & \ts\vb\cr
\vsa
&&&&\multispan{5}{\hrulefill}\cr
\ssa
\vb & 2 & & 8 & \ts\vb\cr
\vsa
&&&&\multispan{0}{\hrulefill}\cr
\ssa
\vb & 4 & & 9 & \ts\vb\cr
\vsa
&&\multispan{3}{\hrulefill}\cr
\ssa
\vb & 5 & \ts\vb\cr
\vsa
&&\multispan{0}{\hrulefill}\cr
\ssa
\vb & 6 & \ts\vb\cr
\vsa
\multispan{3}{\hrulefill}\cr}}\right\},\quad
T_\Iscr(3)=\left\{
\vcenter{
\halign{& \hfill#\hfill
\tabskip4pt\cr
\multispan{9}{\hrulefill}\cr
\ssa
\vb & 1 &  &  3 && 7& & 10 & \ts\vb\cr
\vsa
&&&&&&\multispan{3}{\hrulefill}\cr
\ssa
\vb & 2 &  &  4 & & 8 &\ts\vb\cr
\vsa
&&\multispan{5}{\hrulefill}\cr
\ssa
\vb & 5 &\ts\vb\cr
\vsa
&&\multispan{0}{\hrulefill}\cr
\ssa
\vb & 6 &\ts\vb\cr
\vsa
&&\multispan{0}{\hrulefill}\cr
\ssa
\vb & 9 &\ts\vb\cr
\vsa
\multispan{3}{\hrulefill}\cr}},\quad
\vcenter{
\halign{& \hfill#\hfill
\tabskip4pt\cr
\multispan{9}{\hrulefill}\cr
\ssa
\vb & 1 &  &  3 && 7& & 10 & \ts\vb\cr
\vsa
&&&&\multispan{5}{\hrulefill}\cr
\ssa
\vb & 2 &  &  4 & \ts\vb\cr
\vsa
&&&&\multispan{0}{\hrulefill}\cr
\ssa
\vb & 5 && 8 &\ts\vb\cr
\vsa
&&&&\multispan{0}{\hrulefill}\cr
\ssa
\vb & 6 && 9&\ts\vb\cr
\vsa
\multispan{5}{\hrulefill}\cr}}
\right\}, $$
$$T_\Iscr(4)=\left\{
\vcenter{
\halign{& \hfill#\hfill
\tabskip4pt\cr
\multispan{7}{\hrulefill}\cr
\ssa
\vb & 1 &  &  3 && 7& \ts\vb\cr
\vsa
&&&&&&\multispan{0}{\hrulefill}\cr
\ssa
\vb & 2 &  &  4 & & 8 &\ts\vb\cr
\vsa
&&&&&&\multispan{0}{\hrulefill}\cr
\ssa
\vb & 5 && 9 && 10&\ts\vb\cr
\vsa
&&\multispan{5}{\hrulefill}\cr
\ssa
\vb & 6 &\ts\vb\cr
\vsa
\multispan{3}{\hrulefill}\cr}}
\right\}. $$
{\bf Remark 1.} For $i>1$ one has $T_\Iscr(i)\ne \emptyset.$
\parno 
{\bf Remark 2.} Define
$$\hat T_\Iscr(i,\prev):=\begin{cases}T_\Iscr(i,\prev_l(i),\prev_l(i)), & {\rm if}\ \prev_l(i)\ne n,\cr
                          \emptyset, & {\rm otherwise}.\cr\end{cases}$$
Then
$$\hat T_\Iscr(i,\prev):=\begin{cases}T_\Iscr(i,\next), & {\rm if}\ \prev_l(i)=c_i,\cr
                                T_\Iscr(i,\prev), & {\rm otherwise}.\cr\end{cases}$$
That is why we define $T_\Iscr(i,\prev)$ with the help of $\sprev_l(i)$ and not of $\prev_l(i).$ 
\subsection{}\label{2.6} 
Let $\Dscr_G(T)$ denote the set of geometric descendants of $T.$
Now we can formulate the main theorem of this paper.
\begin{theorem} 
For $T_\Iscr$ the set of descendants is defined by
$$\begin{array}{rl} \Dscr_G(T_\Iscr)=&\{T_\Iscr(i,{\next})\ |\ {\next}_l(i)=c_i\ {\rm or}\ 
\langle{\next}_l(i)>c_i\ {\rm and}\ {\snext}_r(i)>{\next}_l(i)\rangle\}\bigcup\cr
&\{T_\Iscr(i,\prev)\ |\ \sprev_l(i)<n\ {\rm and}\ \sprev_r(i)<\sprev_l(i)\}.\cr\end{array}$$
\end{theorem}
In our example from \ref{2.4} one has
$$\Dscr_G(T_\Iscr)=\left\{
\vcenter{
\halign{& \hfill#\hfill
\tabskip4pt\cr
\multispan{9}{\hrulefill}\cr
\ssa
\vb & 1 &  &  3 && 7& & 10 & \ts\vb\cr
\vsa
&&&&&&\multispan{3}{\hrulefill}\cr
\ssa
\vb & 2 &  &  4 & & 8 &\ts\vb\cr
\vsa
&&\multispan{5}{\hrulefill}\cr
\ssa
\vb & 5 &\ts\vb\cr
\vsa
&&\multispan{0}{\hrulefill}\cr
\ssa
\vb & 6 &\ts\vb\cr
\vsa
&&\multispan{0}{\hrulefill}\cr
\ssa
\vb & 9 &\ts\vb\cr
\vsa
\multispan{3}{\hrulefill}\cr}},\quad
\vcenter{
\halign{& \hfill#\hfill
\tabskip4pt\cr
\multispan{9}{\hrulefill}\cr
\ssa
\vb & 1 &  &  3 && 7& & 10 & \ts\vb\cr
\vsa
&&&&\multispan{5}{\hrulefill}\cr
\ssa
\vb & 2 &  &  4 & \ts\vb\cr
\vsa
&&&&\multispan{0}{\hrulefill}\cr
\ssa
\vb & 5 && 8 &\ts\vb\cr
\vsa
&&&&\multispan{0}{\hrulefill}\cr
\ssa
\vb & 6 && 9&\ts\vb\cr
\vsa
\multispan{5}{\hrulefill}\cr}},
\quad
\vcenter{
\halign{& \hfill#\hfill
\tabskip4pt\cr
\multispan{7}{\hrulefill}\cr
\ssa
\vb & 1 &  &  3 && 7& \ts\vb\cr
\vsa
&&&&&&\multispan{0}{\hrulefill}\cr
\ssa
\vb & 2 &  &  4 & & 8 &\ts\vb\cr
\vsa
&&&&&&\multispan{0}{\hrulefill}\cr
\ssa
\vb & 5 && 9 && 10&\ts\vb\cr
\vsa
&&\multispan{5}{\hrulefill}\cr
\ssa
\vb & 6 &\ts\vb\cr
\vsa
\multispan{3}{\hrulefill}\cr}}
\right\}. $$
Consider  $i:\ 2\leq i\leq l.$ Note that for any 
$T\in T_\Iscr(i)$ one has $\tau(T)\supset \Iscr$ thus
$T\gegs T_\Iscr.$
\par  
The proof of the theorem consists of the following steps:
\begin{itemize}
\item[(i)] We show that
$$\Dscr_G(T_\Iscr)\st \bigcup\limits_{i=2}^l T_\Iscr(i).$$
\item[(ii)] We show that $T\in T_\Iscr(l)$ is always in 
$\Dscr_G(T_\Iscr).$ 
\item[(iii)] Finally we consider $T\in T_\Iscr(i)$ for $i\ :\ 1<i<l.$
Note that if ${\next}_l(i)=c_i$ then $T_\Iscr(i)=\{T_\Iscr(i,{\next})\}$
and $\codim_{\gm_\Iscr} \Vscr_{T_\Iscr(i,\next)}=1,$ thus
$T_\Iscr(i,\next)\in \Dscr_G(T_\Iscr)$ just by the dimension 
consideration. We must consider only $T\in T_\Iscr(i)$
such that $\codim_{\gm_\Iscr} \Vscr_T>1,$ that is
$T_\Iscr(i,\next)$ when $\next_l(i)> c_i$
and $T_\Iscr(i,\prev)$ if it exists. To deal with these cases
we use order preserving involution $\psi$ on $\bT_n$ induced by 
the non-trivial involution of the Dynkin diagram of $\gs\gl_n.$
With the help of $\psi$ we first show
that in case $\next_l(i)> c_i$ one has
$T_\Iscr(i,\next)\in \Dscr_G(T_\Iscr)$ iff $\snext_r(i)>\next_l(i).$
Then we use $\psi$ again to show that $T_\Iscr(i,\prev)\in
\Dscr_G(T_\Iscr)$ iff $\sprev_r(i)<\sprev_l(i).$
\end{itemize}
\section{\bf Proof of the theorem on descendants}
\subsection{}\label{3.0}
For $1\leq i<j\leq n$ set $<i,j>:=\{k\}_{k=i}^j.$
Recall the notion of projection $\pi_{i,j}:{\bT}_n\rar {\boT}_{<i,j>}$ obtained through
jeu de taquin applied to entries of $T\in \bT_n$ not lying in $<i,j>$ as it is defined in I, 2.4.16.
\subsection{}\label{3.1}
Let us consider the projections $\pi_{i,j}$ of $T\in T_\Iscr(i).$
As noted in \ref{2.2} the tableau $\pi_{1,n-1}(T_\Iscr)$
is obtained from $T_\Iscr$ by eliminating $n$ from $C_l.$ It corresponds 
to the subset $\Iscr_n:=\Iscr\cap \Pi_{1,n-1}$
with $C^{\Iscr_n}_i=C^\Iscr_i\ :\ i<l$ and $C^{\Iscr_n}_l=
C^\Iscr_l-n.$ Similarly $\pi_{2,n}(T_{\Iscr})$ is obtained from
$T_\Iscr$ by eliminating 1 from $C_1$ and corresponds to 
the  subset $\Iscr_1:=\Iscr\cap\Pi_{2,n}$ with 
$C^{\Iscr_1}_i=C^\Iscr_i\ :\ i>1$ and $C^{\Iscr_1}_1=C^\Iscr_1-1.$
\par
If $T_\Iscr(i,\next)\ne \emptyset$  set $p(i,\next):=\max\{s\ |\ s<i,\ c_s=\next_l(i)\}.$
If $T_\Iscr(i,\prev)\ne \emptyset$ set $p(i,\prev):=\max\{s\ |\ s<i\ c_s=\sprev_l(i)\}.$
\begin{lemma} For $i: i>1$ consider $T\in T_\Iscr(i)$ and fix 
$$p(T)=\begin{cases} p(i,\next), & {\rm if}\ T=T_\Iscr(i,\next)\cr
              p(i,\prev), & {\rm otherwise}\cr\end{cases}$$
Then
\begin{itemize}
\item[(i)] $\pi_{1,n-1}(T)\in T_{\Iscr_n}(i),$ if $i<l.$
\item[(ii)] If $T_\Iscr(l,\next)\ne \emptyset$ then $\pi_{1,n-1}(T_\Iscr(l,\next))=T_{\Iscr_n}.$
\item[(iii)] $\pi_{1,n-1}(T_\Iscr(l,\prev))\in T_{\Iscr_n}(l).$
\item[(iv)] $\pi_{2,n}(T)\in T_{\Iscr_1}(i),$ if $p(T)>1$
or $p(T)=1$ and $c_1>c_i.$
\item[(v)] $\pi_{2,n}(T)=T_{\Iscr_1},$ if $p(T)=1$ and $c_1\leq c_i.$
\end{itemize}
\end{lemma}
\Pf
(i) is an immediate consequence of $\si_i<n,$ for $i<l.$
\parno
(ii) is an immediate consequence of the definition of $T_\Iscr(l,\next)$
( which is defined by moving the box containing $n$ down) 
and of $\pi_{{\sr 1},n-\sr 1}(T)$ ( which is defined by eliminating box containing $n$ from $T$).
\parno
(iii) is an immediate consequence of $C_l^{\sprev(l)}<n.$
\parno 
To prove (iv) note that $T_{\Iscr_1}$ is obtained by replacing $C_1$ by $C_1-1$ 
If $p(T)>1$ then the operation of moving the tail commutes with  
this replacing. If $p(T)=1$ and $c_1>c_i$ then $T=T_\Iscr(i,\next)$
and $\si_i$ goes from 
$c_i-$th row to  
$(c_1+1)-$th row in forming $T_\Iscr(i,\next),$ but then goes to $c_1-$th row in
forming $\pi_{2,n}(T_\Iscr(i,\next)).$ However if $c_1>c_i$ this is exactly what happens in
forming $T_{\Iscr_1}(i,\next)$ from $T_{\Iscr_1}.$ Hence (iv).
\parno
Finally if $c_i>c_1>c_j,\ \forall j: 1<j<i,$
that is under hypothesis of (v) for $T_\Iscr(i,\prev)$, one has 
$T^{c_1}_2=C_i^{c_1}$
and it is pushed back to $c_1-$th row in forming 
$\pi_{2,n}(T)$ which is hence $T_{\Iscr_1}.$
\QED
\subsection{}\label{3.2}
Recall notation $p_w(i)$ from \ref{2.0}.
In what follows we need the following simple
\begin{lemma} Given $x,z\in \bS_n$ such that
\begin{itemize} 
\item[(i)] $\pi_{1,s}(x)=\pi_{1,s}(z)$  
\item[(ii)] $p_x(n)<p_x(n-1)<\cdots<p_x(s+1),\  p_z(n)<p_z(n-1)<\cdots<p_z(s+1)$
\item[(iii)] $p_x(j)\leq p_z(j)\ \forall j\ : s+1\leq j\leq n.$
\end{itemize}
Then $\gn\cap^x\gn\st \gn\cap^z\gn.$ 
\end{lemma}
\Pf
Recall from \ref{2.0} that $\al_{i,j}\in S(w)$ if and only if
$p_w(i)>p_w(j).$  Applying this proposition to all 3 hypotheses we get
\begin{itemize}
\item[(i)] implies that for any $i,j\ :\ i<j\leq s$ 
one has $\al_{i,j}\in S(x)$ iff $\al_{i,j}\in S(z);$ 
\item[(ii)] implies that for any $i,j\ :\ s+1\leq i<j$ one has $\al_{i,j}\in S(x)$
and $\al_{i,j}\in S(z);$
\item[(iii)] implies that for any $i,j\ :\ i\leq s,\ j\geq s+1$ one has $\al_{i,j}\in S(z)$ implies
$\al_{i,j}\in S(x).$ 
\end{itemize}
Thus $S(z)\st S(x)$ which is equivalent to $\gn\cap^x\gn\st \gn\cap^z\gn.$
\QED
\subsection{}\label{3.3}
Given tableaux $P,Q$ such that $<P>\cap<Q>=\emptyset,$
recall the notation $(P,Q)$ from \ref{2.01} or Part I, 2.4.3 and set
$\binom{P}{Q}:=(P,Q)^{\dagger}.$ 
\par
In what follows we need the lemma which is a kind of variation
of Part I, 3.2.3 (v),(vi).
\begin{lemma} Given tableaux $P$ and $Q$ such that for any 
$p\in <P>$
and any $q\in <Q>$ one has $p<q.$ For any $x\in \Cscr_P$ and 
$y\in \Cscr_Q$
one has 
\begin{itemize}
\item[(i)] Set $w=[x,y].$ Then $T(w)=(P,Q).$
\item[(ii)] Set $w=[y,x].$ Then $T(w)=\binom{P}{Q}.$
\end{itemize}
\end{lemma}
\par
\Pf
To show (i) it is enough to show that for any row $i$ of $P,$ for any
$p\in P^i$
one has $r_{\sr T(w)}(p)=i$ and for any 
row $i$ of $Q,$ for any
$q\in Q^i$
one has $r_{\sr T(w)}(q)=i.$ Let $x=[x_{\sr 1},\ldots,x_{\sr k}]$
and $y=[y_{\sr 1},\ldots,y_m].$
\par
Recall notation of RS procedure from Part I, 2.4.6.
One has ${}_kT(w)=T(x)=P$ so that $r_{\sr {}_kT(w)}(p)=i.$
For any $y_j$  one has $y_j>p$ thus RS 
insertion does not knock down $p$ so that
$r_{\sr T(w)}(p)=i.$ On the other hand 
again since any $y_j$ does not knock down any $x_{j\pr}$
one has that it knocks down exactly the same $Q^i_j$ as in
$T(y).$ Thus if  $q\in Q^i$ then $r_{\sr T(w)}(q)=i.$
\par
To show (ii) we use Schensted-Sch\"utzenberger theorem (cf. Part I, 2.4.15)
claiming $T(\ov w)=(T(w))^{\dagger}.$ First of all this theorem together with
part (i) gives us that $T([\ov x,\ov y])=(P^{\dagger},Q^{\dagger}).$
Applying the theorem again we get 
$$T([y,x])=T(\ov{[\ov x,\ov y]})=(P^{\dagger},Q^{\dagger})^{\dagger}=\binom{P}{Q}.$$
\QED
\subsection{}\label{3.4}
Now we are ready to show that $\Dscr_G(T)\st 
\{T_\Iscr(i)\}_{i=2}^l.$
Let $T\gegs T_\Iscr$ and let
$m$ be minimal integer such that $r_{\sr T}(m)>r_{\sr T_\Iscr}(m).$ 
Let $i$ be the number of the chain $m$ belongs to, that is $m\in<C_i^{\Iscr}>.$
If $T_\Iscr(i)$ contains only one element let $T_\Iscr[m]$ denote this element.
If $T_\Iscr(i)=\{T_\Iscr(i,\next),T_\Iscr(i,\prev)\}$ then 
$$T_{\Iscr}[m]=\begin{cases} T_\Iscr(i,\prev), & {\rm if}\ m\leq C_i^{\sprev_l(i)}\cr
              T_\Iscr(i,\next), & {\rm otherwise}\cr\end{cases}$$
\begin{prop} Let $\ov\Vscr_T\subsetneq\gm_\Iscr$ and let
$m$ be minimal integer such that $r_{\sr T}(m)>r_{\sr T_\Iscr}(m).$ 
Then $\ov\Vscr_T\st\ov\Vscr_{T_\Iscr[m]}.$ In particular  
$\Dscr_G(T)\st \{T_\Iscr(i)\}_{i=2}^l.$
\end{prop}
\Pf
This is trivially true for $\gs\gl_3.$ Assume this is true for
$n-1$ and show this for $n.$
\par
Assume $i<l.$ Consider $\pi_{1,n-1}(T_\Iscr[m]).$ 
By \ref{3.1} (i) in that case
$\pi_{i,n-1}(T_\Iscr[m])=T_{\Iscr_n}[m].$
By Part I, 4.1.2
one has $\pi_{1,n-1}(\ov\Vscr_T)\st \gm_{\Iscr_n}.$ As well,
since $m<n$ in that case, one has $r_{\sr \pi_{1,n-1}(T)}(m)=
r_{\sr T}(m)>r_{\sr T_{\Iscr_n}}(m).$
Thus by induction assumption $\ov\Vscr_{\pi_{1,n-1}(T)}\st
\ov\Vscr_{T_{\Iscr_n}[m]}.$
\par
On the other hand we can write $\gm_\Iscr$ as a Cartesian product: 
$$\gm_{\sr \Iscr}=\gm_{\sr \Iscr_n}\times \sum_{j=1}^{\si_{l-1}}X_{\al_{j,n}}.$$
Thus for any  $\ov\Vscr_T\st \gm_\Iscr$ one has 
$$\ov\Vscr_T\st\ov\Vscr_{\pi_{1,n-1}(T)}\times \sum_{j=1}^{\si_{l-1}}X_{\al_{j,n}}.$$
For any $i: 1<i<l$ consider $w=w_r(\pi_{{\sr 1},\si_{l-1}}(T_\Iscr[m]))$ and 
$x=\ov C_l.$ 
By \ref{3.3} (i) one has
$Q([w,x])=T_\Iscr[m].$ 

By Part I, 2.2.4 one has $\al_{j,n}\in S([w,x])$ if and only if $j\in C_l.$ Thus
$\ov\Vscr_{T_\Iscr[m]}=\ov\Vscr_{T_{\Iscr_n}[m]} \times 
\sum_{j=1}^{\si_{l-1}}X_{\al_{j,n}}.$
Hence if $i<l$ then
$$\ov\Vscr_T\st\ov\Vscr_{\pi_{1,n-1}(T)}\times \sum_{j=1}^{\si_{l-1}}X_{\al_{j,n}}\st 
\ov\Vscr_{T_{\Iscr_n}[m]} \times \sum_{j=1}^{\si_{l-1}}X_{\al_{j,n}}=
\ov\Vscr_{T_\Iscr[m]}$$
\par
If $i=l$ let us show that there exists $y\ :\ T(y)=T$ and 
$z\ :\ T(z)=T_\Iscr[m]$
such that $\gn\cap^y\gn\st\gn\cap^z\gn.$ 
Indeed let $y=w_r(T)$ and $z=w_r(T_\Iscr[m]).$
Note that 
\begin{itemize}
\item[(i)] $\pi_{1,m-1}(T)=\pi_{1,m-1}(T_\Iscr[m])$ therefore 
$\pi_{1,m-1}(y)=\pi_{1,m-1}(z).$
\item[(ii)] $r_{\sr T}(n)>r_{\sr T}(n-1)>\cdots>r_{\sr T}(m)$ and 
$r_{\sr T_\Iscr[m]}(n)>r_{\sr T_\Iscr[m]}(n-1)>\cdots>r_{\sr T_\Iscr[m]}(m)$
therefore $p_y(n)<p_y(n-1)<\cdots<p_y(m)$ and $p_z(n)<p_z(n-1)<\cdots<p_z(m).$
\item[(iii)]\begin{itemize}\item[(a)] If $m\leq C_l^{\sprev_l(l)}$ then for any $s\geq m$ one has
$$r_{\sr T_\Iscr[m]}(s)=
\begin{cases}r_{\sr T_\Iscr}(s),& {\rm if} s<C_l^{\sprev_l(l)},\cr
       r_{\sr T_\Iscr}(s)+1, & {\rm otherwise}.\cr\end{cases}$$
On the other hand for any $s\geq m$ one has $r_{\sr T}(s)\geq r_{\sr T_\Iscr}(s)+1$
just by $\{\al_s\}_{s=m}^{n-1}\st\tau(T_\Iscr)\st\tau(T)$ and by the condition 
$r_{\sr T}(m)\geq r_{\sr T_\Iscr}(m)+1.$  
Thus $p_y(s)\leq p_z(s).$
\item[(b)] If $m> C_l^{\sprev_l(l)}$ then 
$$r_{\sr T_\Iscr[m]}(s)=\begin{cases}r_{\sr T_\Iscr}(s),& {\rm if}\ s<n,\cr
                   \snext_l(l)+1,& {\rm if}\ s=n.\cr\end{cases}$$
On the other hand since $|T_\Iscr^{\sprev_l(l)+1}|=\cdots=|T_\Iscr^{c_l}|$
and $|T_\Iscr^{c_l+1}|=\cdots=|T_\Iscr^{\snext_l(l)}|=|T_\Iscr^{c_l}|-1$
one has that $r_{\sr T}(m)\geq \snext_l(l)+1$ and 
for any $s\geq m$ one has $r_{\sr T}(s)\geq \snext_l(l)+1+(s-m)$ exactly by the same
reasoning as in (a). Thus again $p_y(s)\leq p_z(s).$
\end{itemize}
\end{itemize}
In both cases we get the hypothesis (iii) of \ref{3.2}.
Therefore by \ref{3.2} we get the result.
\QED
\parno  
{\bf Remark.}\ \ Recall that $T_\Iscr[m]\in T_\Iscr(i).$ 
Let $w=w_r(\pi_{1,\si_i}(T_\Iscr[m]))$
and 
$$x=\begin{cases}[\ov C_{i+1},\ldots,\ov C_l]& {\rm if}\ i<l,\cr
            \emptyset, & {\rm if}\ i=l.\cr\end{cases}$$ 
then by \ref{3.3} (i) one has $Q[w,x]=T_\Iscr[m]$ and the same proof shows that 
$$\gn\cap^{w_r(T)}\gn\st \gn\cap^{[w,x]}\gn,$$ 
that is $T\dgs T_\Iscr[m].$
In particular this  means that $T\gegs T_\Iscr[m]$ 
if and only if 
$T\dgs T_\Iscr[m].$ 
The interesting question is whether this is true
for any $T\in\{T_\Iscr(i)\}_{i=2}^l,$ that is whether for any $T\in\{T_\Iscr(i)\}_{i=2}^l$
one has $S\gegs T$ if and only if $S\dgs T.$ 
\subsection{}\label{3.5}
We have completed step (i) of the proof. Now we prove step (ii), namely 
\begin{prop} Given $T_\Iscr=(C_1,\ldots,C_l).$
\begin{itemize}
\item[(i)] If $T_\Iscr(l,\next)\ne \emptyset$ then $T_\Iscr(l,\next)\in \Dscr(T_\Iscr).$
\item[(ii)] If $T_\Iscr(l,\prev)\ne\emptyset$ then $T_\Iscr(l,\prev)\in \Dscr(T_\Iscr).$
\end{itemize}
\end{prop}
\Pf
Set $\sh(T_\Iscr)=\lam=(\lam_1,\lam_2,\ldots,\lam_k).$
\par 
(i)\ \ \  Assume that $T_\Iscr(l,\next)\ne \emptyset$
and put $T:=T_\Iscr(l,\next),\ l\pr:=p(l,\next)$ and $\sh(T):=\mu=(\mu_1,\ldots,\mu_k).$ 
(If the number of rows in $T$ is greater by 1 than the number of rows in $T_\Iscr$ we suppose
that indeed $\lam=(\lam_1,\lam_2,\ldots,\lam_{k-1})$ and add 
$\lam_k=0$ to $\lam$ to get the same length.) 
One has
$$\mu_j=\begin{cases}\lam_j, & {\rm if}\ j\ne c_l,c_{l\pr}+1,\cr
               \lam_{c_l}-1, & {\rm if}\ j=c_l,\cr
               \lam_{c_{l\pr}+1}+1, & {\rm if}\ j=c_{l\pr}+1.\cr\end{cases}$$
If $c_{l\pr}=c_l$ then $\Oscr_\mu$ is a descendant
of $\Oscr_\lam$ by Part I, 2.3.3. Thus $\Vscr_T$ is a descendant
of $\Vscr_\Iscr.$ If for any $j\ne l\pr$ one has $c_j\ne c_{l\pr}$
then $\lam_{c_l+1}=\cdots=\lam_{c_{l\pr}}=\lam_{c_l}-1$ and $\lam_{c_{l\pr}+1}=\lam_{c_l}-2$
so that again $\Oscr_\mu$ is a descendant
of $\Oscr_\lam$ by Part I, 2.3.3. Thus again $\Vscr_T$ is a descendant
of $\Vscr_\Iscr.$
\par
Assume that $c_{l\pr}>c_l$ and there exist $j\ne l\pr$ such that $c_j=c_{l\pr}.$
Then there exist the unique orbit $\Oscr_\nu$ such that 
$\ov\Oscr_\mu\subsetneq \ov\Oscr_\nu\subsetneq\ov\Oscr_\lam$ where $\nu=(\nu_1,\ldots,\nu_k)$
is defined by
$$\nu_j=\begin{cases}\lam_j, & {\rm if}\ j\ne c_{l\pr},c_{l\pr}+1,\cr
               \lam_{c_{l\pr}}-1, &  {\rm if}\ j=l\pr,\cr
               \lam_{c_{l\pr}+1}+1, &  {\rm if}\ j=c_{l\pr}+1.\cr\end{cases}$$
Assume that there exists $T\pr$ such that $\ov\Vscr_T\subsetneq \ov\Vscr_{T\pr}\subsetneq \ov\Vscr_\Iscr.$
Then $\sh(T\pr)=\nu.$ Consider $\pi_{1,n-1}$ of the three orbital variety closures. We get 
$\ov\Vscr_{\pi_{1,n-1}(T)}\st \ov\Vscr_{\pi_{1,n-1}(T\pr)}\st \ov\Vscr_{\pi_{1,n-1}(\Iscr)}.$
Since $\pi_{1,n-1}(T)=\pi_{1,n-1}(T_\Iscr)$ one has $\pi_{1,n-1}(T\pr)=\pi_{1,n-1}(T_\Iscr)$
which is impossible since
$|(\pi_{1,n-1}(T_\Iscr))^{c_{l\pr}}|=\lam_{c_{l\pr}}$ and 
$|(\pi_{1,n-1}(T\pr))^{c_{l\pr}}|\leq \nu_{c_{l\pr}}=\lam_{c_{l\pr}}-1.$
\par
(ii)\ \ \ The proof of (ii) is very similar to the proof
of (i). Assume that $T_\Iscr(l,\prev)\ne\emptyset$ and put $T=T_\Iscr(l,\prev),$
$l\pr:=p(l,\prev)$ and $\sh(T):=\mu=(\mu_1,\ldots,\mu_k).$ 
(Again if the number of rows in $T$ greater by 1 than the number of rows in $T_\Iscr$  we again
add $\lam_k=0$ to $\lam$ to get the same length.) One has
$$\mu_j=\begin{cases}\lam_j, &  {\rm if}\ j\ne c_{l\pr},c_l+1,\cr
               \lam_{c_{l\pr}}-1, &  {\rm if}\ j=c_{l\pr},\cr
               \lam_{c_l+1}+1, &  {\rm if}\ j=c_l+1.\cr\end{cases}$$ 
Note that by definition of $T_\Iscr(l,\prev)$ we obtain 
$$\lam_{c_{l\pr}+1}=\cdots=\lam_{c_l}\quad {\rm and}\quad \lam_{c_l+1}=\lam_{c_l}-1.\eqno{(*)}$$
The claim is trivially true for $n=3$ so assume that in $\bT_{n-1}$ one has
that $T_\Iscr(l,\prev)$ is a descendant of $T_\Iscr.$
\par
Again, if for any $j\ne l\pr$ one has $c_j\ne c_{l\pr}$
then $\Oscr_\mu$ is a descendant of $\Oscr_\lam$ by Part I, 2.3.3, thus $\Vscr_T$ is a descendant
of $\Vscr_\Iscr.$ 
\par
Again, assume that $c_{l\pr}<c_l$ and there exist $j\ne l\pr$ such that $c_j=c_{l\pr}.$
Note that in that case in particular
$$\al_{n-1}\in \Iscr \eqno{(**)}.$$
By $(*)$ in that case there exist the unique orbit $\Oscr_\nu$ such that 
$\ov\Oscr_\mu\subsetneq \ov\Oscr_\nu\subsetneq\ov\Oscr_\lam$ where $\nu=(\nu_1,\ldots,\nu_k)$
is defined by
$$\nu_j=\begin{cases}\lam_j, &  {\rm if}\ j\ne c_{l\pr},c_{l\pr}+1,\cr
               \lam_{c_{l\pr}}-1, &  {\rm if}\ j=l\pr,\cr
               \lam_{c_{l\pr}+1}+1, &  {\rm if}\ j=c_{l\pr}+1.\cr\end{cases}$$
Assume that there exist $T\pr$ such that 
$\ov\Vscr_T\subsetneq \ov \Vscr_{T\pr}
\subsetneq \ov \Vscr_\Iscr$ then $\sh(T\pr)=\nu.$  Consider $\pi_{1,n-1}$ of the three orbital 
variety closures. We get 
$\ov\Vscr_{\pi_{1,n-1}(T)}\st \ov\Vscr_{\pi_{1,n-1}(T\pr)}\st \ov\Vscr_{\Iscr_n}.$
By the induction assumption $\pi_{1,n-1}(T\pr)=\pi_{1,n-1}(T)$ or $\pi_{1,n-1}(T\pr)=T_{\Iscr_n}.$
The second situation is impossible since $|T_{\Iscr_n}^{c_{l\pr}}|=\lam_{c_{l\pr}}$ and 
$|(\pi_{1,n-1}(T\pr))^{c_{l\pr}}|\leq \nu_{c_{l\pr}}=\lam_{c_{l\pr}}-1.$ Thus 
$\pi_{1,n-1}(T\pr)=\pi_{1,n-1}(T).$ Set $\mu\pr=\sh(\pi_{1,n-1}(T))=\sh(\pi_{1,n-1}(T\pr)).$ Then
$$\mu\pr_j=\begin{cases}\lam_j, &  {\rm if}\ j\ne c_{l\pr},\cr
                \lam_{c_{l\pr}}-1, &  {\rm if}\ j=c_{l\pr}.\cr
\end{cases}=\begin{cases}\nu_j, &  {\rm if}\ j\ne c_{l\pr}+1 ,\cr
                         \nu_{c_{l\pr}+1}-1, & {\rm if}\ j=c_{l\pr}+1.\cr\end{cases}$$
Thus $r_{\sr T\pr}(n)=c_{l\pr}+1\leq r_T(n-1)=c_l.$ Since $\pi_{1,n-1}(T\pr)=\pi_{1,n-1}(T)$ we get
that $r_{\sr T}(n-1)=r_{\sr T\pr}(n-1)$ so that $\al_{n-\sr 1}\not\in\tau(T\pr)$ 
which contradicts to $(**).$ 
\QED
\subsection{}\label{3.6}
Note that \ref{3.5} (i) can be easily generalized to any Young tableau $T$
\begin{prop} Given $T\in \bT_n.$ Let $\sh(T)=\lam=(\lam_1,\ldots\lam_k).$ 
Let $i:=r_{\sr T}(n)$ and let $j=\min\{m\ |\ \lam_m<\lam_i-1\}.$ Then $S$ obtained
from $T$ by moving box with $n$ from row $i$ to row $j,$ or formally
$$S:=\left( \begin{array}{c}T^{1,i-1}\cr
           T^i-n\cr
           T^{i+1,j-1}\cr
           T^j+n\cr
           T^{j+1,\infty}\cr\end{array}\right),$$
is a geometric descendant of $T.$
\end{prop}
The proof is exactly the same as the proof of \ref{3.5} (i).
In particular \ref{3.5} and \ref{3.6} show that $\Vscr_S$ being a geometric descendant of $\Vscr_T$
does not necessarily imply that $\Oscr_S$ is a descendant of $\Oscr_T.$
\subsection{}\label{3.7}
As it is formulated in theorem \ref{2.6}, it may happen that for a 
given $I\ :\ i<l,\ \ T_\Iscr(i,\next)$ or 
$T_\Iscr(i,\prev)$ is {\it not} a geometric
{\it descendant} of $T_\Iscr.$ This is a difficult point. 
To understand it we need the non-trivial involution 
$\psi$ of  the Dynkin diagram of $\gs\gl_n$ defined by 
$\psi(\al_i)=\al_{n-i}.$ This involution induces the involution
$\psi$ of $\gs\gl_n$ obtained by $\psi(X_{\al_{i,j}})=X_{\al_{n+1-j,n+1-i}}$ 
and of its Weyl group $W=\bS_n$ 
defined by $\psi(s_{\al_i})=s_{\al_{n-i}}.$ One has $\psi(w)(\psi(\al_{i,j}))=\psi(w(\al_{i,j})).$
The involution has a nice description
on word presentations. Set $\psi_n:\{i\}_{i=1}^n\rar \{i\}_{i=1}^n$ by
$\psi_n(i):=n+1-i.$
\begin{lemma} Let $w=[a_1,a_2,\ldots a_n]$ be a word presentation.

Then $\psi(w)=[\psi_n(a_n),\psi_n(a_{n-1}),\ldots,\psi_n(a_1)].$
\end{lemma}
\Pf
First of all let us show the assertion for $s_i=[1,2,\ldots, i-1,i+1,i,i+2,\ldots,n]$
$$\begin{array}{rl}
\psi(s_i)=s_{n-i}&=[1,2,\ldots,n-i-1,n-i+1,n-i,\ldots,n]\cr
&=
[\psi_n(n),\ldots,\psi_n(i+2),\psi_n(i),\psi_n(i+1),\psi_n(i-1),\ldots,\psi_n(1)].\cr
\end{array}$$
Recall from Part I, 1.9 the notion of $\ell(w)$ -- that is the minimal length of $w$ as an element
of Weyl group written as a product of $s_{\al}\ :\ \al\in\Pi.$
Assume that the assertion is true for any $w$ of length
$\ell(w)\leq p-1$ and show this for 
\break
$w=s_{i_1}\cdots s_{i_p}.$ Let $w\pr=s_{i_1}\cdots s_{i_{p-1}}=
[a_1,\ldots a_n]$ in word presentation. Then 
\break
$w=w\pr s_{i_p}=[b_1,\ldots,b_n]$ where
$$b_j=\begin{cases}a_j, & {\rm if}\ j\ne i_p,i_p+1,\cr
             a_{i_p+1}, & {\rm if}\ j=i_p,\cr
             a_{i_p}, & {\rm if}\ j=i_p+1.\cr\end{cases}$$
One has
$$\begin{array}{rlr}\psi(w)&=\psi(w\pr)s_{n-i_p}&\cr
&{=[\psi_n(a_n),\ldots,\psi_n(a_{i_p}),\psi_n(a_{i_p+1}),\ldots,\psi_n(a_1)]s_{n-i_p}}&{
{\rm by\ ind. \ assumption}}\cr
&{=[\psi_n(a_n),\psi_n(a_{n-1}),\ldots,\psi_n(a_{i_p+1}),\psi_n(a_{i_p}),\ldots,\psi_n(a_1)]}&\cr
&{=[\psi_n(b_n),\psi_n(b_{n-1}),\ldots,\psi_n(b_{i_p}),\psi_n(b_{i_p+1}),\ldots,\psi_n(b_1)]}& \cr
\end{array}$$
\QED
\par
In what follows we use the following notation. For
$w=[a_{\sr 1},\ldots,a_m]$ such that
$<w>\st\{i\}_{i=1}^n$ put $\psi(w):=[\psi_n(a_m),\psi_n(a_{m-\sr 1}),\ldots,\psi_n(a_{\sr 1})].$
\subsection{}\label{3.8}
The involution $\psi$ on $W$  induces the (geometric and Duflo) order preserving 
involution on $\bT_n$
obtained by $\psi(T(w)):=T(\psi(w)).$ Indeed since this is an involution on $\gs\gl_n$ and
on $\bB$ and  
$\psi(\gn\cap^w\gn)=\gn\cap^{\psi(w)}\gn$
we obtain
$$\gn\cap^w\gn\subseteq \ov{\bB(\gn\cap^y\gn)} \quad {\rm iff}\quad 
\gn\cap^{\psi(w)}\gn\subseteq \ov{\bB(\gn\cap^{\psi(y)}\gn)}.$$
\par
In general there
is no simple straightforward description of $\psi(T)$ for a given $T.$
But for $T=(P,Q)$ such that $\forall p\in<P>,\ \forall q\in<Q>$ one has $p<q,$ the result
is very simple. 
\begin{lemma} Given tableaux $P$ and $Q$ such that $\forall p\in<P>,\ \forall q\in<Q>$ one has $p<q,$ then
$\psi(P,Q)=(\psi(Q),\psi(P)).$ In particular
for $T_{\Iscr}=(C_1,\ldots,C_l)$ 
one has $\psi(T_{\Iscr})=(\psi(C_l),\ldots,\psi(C_1)).$  
\end{lemma}
\Pf
This is straightforward from \ref{3.3} and \ref{3.7}. 
Let $P=T(x),\ Q=T(y).$ Note that for any $a\in \psi(y)$ and any $b\in \psi(x)$ one has $a<b.$
Now by \ref{3.3} $(P,Q)=T([x,y])$ and 
$$\begin{array}{rlr}\psi(P,Q)&{=T(\psi([x,y]))=T([\psi(y),\psi(x)])}&
{\rm by\ \ref{3.7}}\cr
&{=(\psi(Q),\psi(P))}&{\quad {\rm by\ \ref{3.3}.}}\cr\end{array}$$
\QED
Note that $\psi(T_{\Iscr})=(\psi(C_l),\ldots,\psi(C_1))=T_{\psi(\Iscr)}$ can be also obtained
as a straightforward corollary of the involution $\psi$ of Dynkin diagram.
\subsection{}\label{3.9}
As well we can describe explicitly $\psi(T_\Iscr(i,\next))$ and $\psi(T_\Iscr(i,\prev)).$
This requires more subtle consideration of RS insertion. 
Recall the notation $p(i,\next)$ and $p(i,\prev)$ from \ref{3.1}. Recall the definition
$T_\Iscr(c_s\diar C_i(j))$ from \ref{2.4}.
\begin{prop} Let $T_\Iscr=(C_1,C_2,\ldots,C_l).$ 
\begin{itemize}
\item[(i)] If $T_\Iscr(i,\next)\ne\emptyset$ 
then $\psi(T_\Iscr(i,\next))=T_{\psi(\Iscr)}(c_i\diar \psi(C_{p(i,\next)})(c_i)).$
\item[(ii)] If $T_\Iscr(i,\prev)\ne\emptyset$ then 
$\psi(T_\Iscr(i,\prev))=T_{\psi(\Iscr)}(c_i\diar \psi(C_{p(i,\prev)})(c_{p(i,\prev)})).$
\end{itemize}
\end{prop}
\Pf
Let us first show (i). 
\parno
Set $T:=T_\Iscr(i,\next).$ Note that 
$\sh(T)=\sh(T_{\psi(\Iscr)}(c_i\diar \psi(C_{p(i,\next)})(c_i))).$
\par
First assume that $i=l$ and $p(l,\next)=1.$ In particular $\next_l(l)=c_{\sr 1}.$ 
Note that for any $\Iscr$ for any $i$ one has $|T^i|$ is the number of chains of length
greater or equal than $i.$
In particular by our assumption for $j>1$ if $c_j>c_l$ then $c_j>c_1.$ Thus
$|T_\Iscr^{c_l+1}|=|T^{c_l+2}_\Iscr|=\cdots|T^{c_1}_\Iscr|=:s$ and
$|T_\Iscr^{c_l}|=s+1,\ |T_\Iscr^{c_1+1}|=s-1.$ Recall that $T$ is obtained from
$T_\Iscr$ by moving $n$ from row $c_l$ to row $c_{\sr 1}+1$ so that
$|T^{c_l}|=|T^{c_l+1}|=\cdots=|T^{c_1}|=|T^{c_1+1}|=s.$
Let $\{C_1,C_{i_2},\cdots C_{i_s}\}$ be the chains of length
greater than $c_l.$
\par
Just for the simplicity of further notation note that $C_1=(1,2,\ldots,c_1),$
 $\psi(C_1)=(n+1-c_1,\ldots,n)$ and $(\psi(C_1))^{c_l+t}=n-c_1+c_l+t$ for any $t\ :\
0\leq t\leq c_1-c_l$
\par
Recall $C_i^{1,j}$ from  \ref{2.4}. If 
$c_i\leq j$ then let $C_i^{1,j}=C_i.$ 
\parno
Set $x=[\ov{C_1^{\mathstrut 1,c_1}},\ov{C_2^{\mathstrut 1,c_1}},\ldots,
\ov{C_{l-1}^{\mathstrut 1,c_1}},\ov{C^{\mathstrut}_l-n}]$, $y=[T^{c_1+1}]$ and $z=w_r(T^{c_1+2,\infty}).$ 
Thus by Part I, 3.2.3 (v) one has $T([z,y,x])=T.$
\par
Now consider $S:=\psi(T)=T(\psi([z,y,x])).$ By \ref{3.7} one has
$\psi([y,x,w])=[\psi(x),\psi(y),\psi(z)].$ 
\begin{itemize}
\item[(1)] Put $S\pr=T([\psi(x)])$ then by \ref{3.8} one has
$S\pr=(\psi(C_l-n),\psi(C_{l-1}^{1,c_1}),\ldots,\psi(C_1^{1,c_1})).$
and in particular $\{(n+1-T^{c_1}_s)<(n+1-T^{c_1}_{s-1})<\cdots<(n+1-T^{c_1}_1)\}\st <(S\pr)^1>$
and $(S\pr)^{c_l}=((n+1-C^{c_1+1-c_l}_{i_s}),\ldots,(n+1-C^{c_1+1-c_l}_{i_2}),n-c_1+c_l).$
\item[(2)] $\psi(y)=[1,n+1-T^{c_1+1}_{s-1},\ldots,n+1-T^{c_1+1}_1]=:[y_{\sr 1},\ldots,y_s].$ Note that
$y_{\sr 1}<y_{\sr 2}<\cdots<y_s$ and $y_j<(n+1-T^{c_1}_{s+1-j})$ for any $j\ :\ 1\leq j\leq s$  
\item[(3)] Set $S\prpr=T([\psi(x),\psi(y)]).$ Then by RS procedure and by (1) and (2) one has
$(S\prpr)^{c_l+1}=(S\pr)^{c_l}.$ In particular $r_{\sr S\prpr}(n-c_1+c_l)=c_l+1.$ 
\end{itemize}
By RS procedure for any $p\in S\prpr$ one has $r_{\sr S}(p)\geq r_{\sr S\prpr}(p).$ 
\parno
In particular
$r_{\sr S}(n-c_1+c_l)\geq c_l+1.$ Since $S\gegs T_{\psi(\Iscr)}$ we get just by $\tau$-invariant
that 
$$r_{\sr S}(n-c_1+c_l+t)\geq c_l+1+t\quad {\rm for\ any}\quad  t\ :\ 0\leq t\leq c_1-c_l.\eqno{(*)}$$ 
As well one has $\sh(S)=\sh(T).$ The only tableau greater (in geometric order) than $T_{\psi(\Iscr)}$
fulfilling $(*)$ and the shape condition 
is 
$$S=T_{\psi(\Iscr)}(c_l\diar \psi(C_1)(c_l)).\eqno{(**)}$$
\par
Now assume that $i<l$ or $j:=p(i,\next)>1.$ Then $T=(P,S,Q)$ where
$$P=(C_1,\ldots, C_{j-1}),\quad S=T_{\pi_{\si_{j-1}+1,\si_i}(\Iscr)}(i-j,\next),\quad Q=(C_{i+1},\ldots,C_l).$$
By \ref{3.8} we get $\psi(T)=(\psi(Q),\psi(S),\psi(P)).$ By (i) we get that 
$$\psi(S)=T_{\psi(\pi_{\si_{j-1}+1,\si_i}(\Iscr))}(c_i\diar \psi(C_j)(c_i)).$$
Thus $\psi(T)=T_{\psi(\Iscr)}(c_i\diar \psi(C_j)(c_i)).$
\bigskip
To show (ii) we set $T:=T_\Iscr(i,\prev)$ and again begin with  the case $i=l$
and $p(l,\prev)=1.$ In that case $T=T_\Iscr(c_1\diar C_l(c_1)).$ Note that this case is dual to (i), that is
$T$ is $S$ from $(**).$ Thus since $\psi$ is involution we get
$$\psi(T)=T_{\psi(\Iscr)}(l,\next)=T_{\psi(\Iscr)}(c_l\diar \psi(C_1)(c_1)).$$ 
We proceed as in part (ii) to obtain the result for any $i$ and $p(i,\prev).$
\QED
\subsection{}\label{3.10}
As we have mentioned already in \ref{2.6} if $\next_l(i)= c_i$
then $\sh (T_\Iscr(i,\next))$ is obtained from $\sh(T_\Iscr)$ by moving a box
from row $c_i$ to row $c_i+1$ so that $\codim_{\gm_\Iscr}(\Vscr_{T_\Iscr(i,\next)})=1.$
Thus in this case $T_\Iscr(i,\next)$ is a (geometric) descendant of $T_\Iscr$ just
by dimension consideration.
\par
To complete the proof of the theorem it remains to show
\begin{prop} \begin{itemize} 
\item[(i)] If $\next_l(i)> c_i$ then
$T_\Iscr(i,\next)\in \Dscr(T_\Iscr)$ iff $\snext_r(i)>\next_l(i).$
\item[(ii)]  If $\prev_l(i)=\sprev_l(i)\ne n$ then
$T_\Iscr(i,\prev)\in \Dscr(T_\Iscr)$ iff $\sprev_r(i)<\sprev_l(i).$
\end{itemize}
\end{prop}
\Pf
We begin with (i). If $T_\Iscr(i,\next)$ is not a descendant, 
that is if there exists $T\ :\ T_\Iscr\gegs T\gegs T_\Iscr(i,\next),$
then by \ref{3.4} $T\in T_\Iscr(j)$ for some $j.$ Moreover by \ref{3.5} (i) $j>i.$
Indeed consider 
$\pi_{1,C_i^{c_i}}(T),\ \pi_{1,C_i^{c_i}}(T_\Iscr(i,\next)).$
Since $\pi_{1,C_i^{c_i}}(T)\ne\pi_{1,C_i^{c_i}}(T_\Iscr(i,\next))$
one has by \ref{3.5} (i) that $\pi_{1,C_i^{c_i}}(T)=\pi_{1,C_i^{c_i}}(T_\Iscr).$
Thus $j>i.$ 
\par
Set $\sh(T_\Iscr):=\lam=(\lam_1,\ldots,\lam_k)$ then $\sh(T_\Iscr(i,\next))=\mu$ where
$$\mu_s=\begin{cases}\lam_s, & {\rm if}\ s\ne c_i,c_{p(i,\next)}+1,\cr
               \lam_{c_i}-1, & {\rm if}\ s=c_i,\cr
               \lam_{c_{p(i,\next)}+1}+1, & {\rm if}\ s=c_{p(i,\next)}+1.\cr\end{cases}$$
Let us first show that if $\snext_r(i)>\next_l(i)$ then $T_\Iscr(i,\next)\in \Dscr(T_\Iscr).$ 
Assume that $\snext_r(i)>\next_l(i).$ That means
that $\forall j:\ j>i$ one has $c_j\leq c_i$ or 
$c_j>c_{p(i,\next)}.$ Then for any $j>i$ if $T_\Iscr(j,\next)$ is defined
then $\sh(T_\Iscr(j,\next))=\nu$ where
$$\nu_s=\begin{cases}\lam_s, & {\rm if}\ s\ne c_j,c_{p(j,\next)}+1;\cr
               \lam_{c_j}-1, & {\rm if}\ s=c_j;\cr
               \lam_{c_{p(j,\next)}+1}+1, & {\rm if}\ s=c_{p(j,\next)}+1.\cr\end{cases}$$
\begin{itemize}
\item[(1)] If $c_j<c_i$ or $c_j>c_{p(i,\next)}$ then $\nu\not<\mu.$ 
\item[(2)] If $c_j=c_i$ put $s=C^1_i,\ t=\si_j$ and consider $\pi_{s,t}.$ Using consequently \ref{3.1} we get
$$\pi_{s,t}(T_\Iscr(i,\next))=T_{\pi_{s,t}(\Iscr)},
\quad \pi_{s,t}(T_\Iscr(j,\next))=T_{\pi_{s,t}(\Iscr)}(j-i+1,\next)$$
\end{itemize}
Hence in both cases $T_\Iscr(i,\next)\not\gegs T_\Iscr(j,\next).$
\par
Again, for any $j>i$ if $T_\Iscr(j,\prev)$ is defined
then $\sh(T_\Iscr(j,\prev))=\nu$ where
$$\nu_s=\begin{cases}\lam_s, & {\rm if}\ s\ne c_j+1,c_{p(j,\prev)};\cr
               \lam_{c_{p(j,\prev)}}-1, & {\rm if}\ s=c_{p(j,\prev)};\cr
               \lam_{c_j+1}+1, & {\rm if}\ s=c_j+1.\cr\end{cases}$$
and we are left only with the cases when $c_j\ne c_i$ (if $c_j=c_i$ then
$T_\Iscr(j,\prev)=\emptyset$) so that $\nu\not<\mu.$
\par
To complete (i) we must show that if $\snext_r(i)\leq \next_l(i)$
then $T_\Iscr(i,\next)\not\in \Dscr(T_\Iscr).$ 
We use the involution $\psi.$
Consider $T_{\psi(\Iscr)}.$ Since $\psi$ is order preserving one has
that $T_\Iscr(i,\next)\in \Dscr(T_\Iscr)$ iff  $\psi(T_\Iscr(i,\next))\in \Dscr(T_{\psi(\Iscr)}).$
Let us show that the last assertion is not true in our case.
Indeed, by \ref{3.9} 
$$\psi(T_\Iscr(i,\next))=T_{\psi(\Iscr)}(c_i\diar \psi(C_{p(i,\next)})(c_i)).$$ 
By \ref{3.4} 
$$T_{\psi(\Iscr)}(c_i\diar \psi(C_{p(i,\next)})(c_i))\geg T_{\psi(\Iscr)}[m]$$ 
where
$m=\psi_{n+1}(C_{p(i,\next)}^{c_i}).$
Note the the number of the chain $\psi(C_{p(i,\next)})$ in $T_{\psi(\Iscr)}$ is $l+1-p(i,\next).$
Since $\snext_r(i)\leq \next_l(i)$ we get that $\prev_l(l+1-p(i,\next))>c_i$ in $T_{\psi(\Iscr)}.$
One has $T_{\psi(\Iscr)}[m]=T_{\psi(\Iscr)}(l+1-p(i,\next),\prev)\ne T_{\psi(\Iscr)}(c_i\diar \psi(C_{p(i,\next)})(c_i)).$
Thus by \ref{3.4} $T_{\psi(\Iscr)}(c_i\diar \psi(C_{p(i,\next)})(c_i))\gegs T_{\psi(\Iscr)}[m]$ which implies 
$T_{\psi(\Iscr)}(l+1-p(i,\next),c_i,c_i)\not\in\Dscr(T_{\psi(\Iscr)}).$
\bigskip
(ii) can be obtained in the same manner as (i) or by applying $\psi$ to (i).
We will show the second way which shows as well that $T_\Iscr(i,\next)\in\Dscr(T_\Iscr)$ if and only if
$T_{\psi(\Iscr)}(l+1-p(i,\next),\prev)\in \Dscr(T_{\psi(\Iscr)}).$ 
\par
Put $\psi(T_\Iscr):=(C\pr_1,\ldots,C\pr_l).$ Note that $\psi(C_s)=C\pr_{l+1-s}$
and $|C\pr_{l+1-s}|=c_s$
for any $s\ :\ 1\leq s\leq l.$ Put $s\pr:=l+1-s$ so that $\psi(C_s)=C\pr_{s\pr}.$
Now set
$$\begin{array}{rl}\snext_l(s)&:=\begin{cases}
\min\{c_j\ |\ j<s,\ c_j> c_s\},& {\rm if\ there\ exists}\ j<s\ {\rm such\ that}\ c_j\geq c_s;\cr
                     n, & {\rm otherwise}.\cr\end{cases}\cr  
  \next_r(s)&:=\begin{cases}\min\{c_j\ |\ j>s,\ c_j\geq c_s\},& {\rm if\ there\ exists}\ j>s 
{\rm such\ that}\ c_j> c_s;\cr
                     0, & {\rm otherwise}.\cr\end{cases}\cr\end{array}$$
Note that for any $s\ : 1\leq s\leq l$
\begin{itemize}
\item[(a)] $\snext_r(s\pr)|_{\psi(\Iscr)}=\snext_l(s)|_\Iscr$ and $\next_l(s\pr)|_{\psi(\Iscr)}=\next_r(s)|_\Iscr.$
\item[(b)] Assume that $s$ is such that $\next_l(s\pr)|_{\psi(\Iscr)}\ne c_s,0.$ 
Set $p\pr:=p(s\pr,\next)|_{\psi(\Iscr)}.$ Respectively we get $p=l+1-p\pr.$ 
\begin{itemize}
\item[(b1)] If $\sprev_l(p\pr)|_{\psi(\Iscr)}\ne n$ then $\sprev_l(p\pr)|_{\psi(\Iscr)}=\sprev_r(p)|_\Iscr$ and
$\sprev_l(p\pr)|_{\psi(\Iscr)}<c_s$ (otherwise $\next_l(s\pr)|_{\psi(\Iscr)}<c_p$). If $\sprev_l(p\pr)|_{\psi(\Iscr)}=n$
then $\sprev_r(p)|_\Iscr=0.$ In both cases we get $\sprev_r(p)|_\Iscr<c_s.$
\item[(b2)] One has $\prev_l(p)|_\Iscr=\sprev_l(p)|_\Iscr=c_s$ iff $\snext_l(s)|_\Iscr>\next_r(s)|_\Iscr.$
\end{itemize} 
\end{itemize}
Consider $T_\Iscr(i,\prev)\ne\emptyset.$ Set $j=p(i,\prev).$
\par  
First assume that $T_\Iscr(i,\prev)\in \Dscr(T_\Iscr).$ In that case 
one has $\psi(T_\Iscr(i,\prev))\in \Dscr(T_{\psi(\Iscr)}).$ By \ref{3.9} (ii)
one has $\psi(T_\Iscr(i,\prev))=T_{\psi(\Iscr)}(c_i\diar C\pr_{j\pr}(c_j)).$
By our assumption one necessarily has by \ref{3.4} that
$\psi(T_\Iscr(i,\prev))=T_{\psi(\Iscr)}(C\pr_{j\pr},\next)$ and that $\next_l(j\pr)|_{\psi(\Iscr)}=c_i>c_j.$
Hence the assumptions of (b) are satisfied for $s=j.$ 
Note that $p(j\pr,\next)_{\psi(\Iscr)}=i\pr$ (otherwise $\prev_l(i)|_\Iscr\ne c_j$).
Using (b1) to $s=j$ and $p=i$ we get that 
$$\sprev_r(i)|_\Iscr<c_j.\eqno{(*)}$$
By (i) one has $T_{\psi(\Iscr)}(C\pr_{j\pr},\next)\in\Dscr(T_{\psi(\Iscr)})$ iff
$\snext_r(j\pr)|_{\psi(\Iscr)}>\next_l(j\pr)|_{\psi(\Iscr)}.$  By (a) this is equivalent to 
$\snext_l(j)|_\Iscr>\next_r(j)|_\Iscr.$ By (b2) this provides $\sprev_l(i)|_\Iscr=c_j.$ 
Comparing this to $(*)$ we get $\sprev_r(i)|_\Iscr<\sprev_l(i)|_\Iscr.$
\par
Now assume that $\sprev_r(i)|_\Iscr<\sprev_l(i)|_\Iscr.$ Recall that $j=p(i,\prev)|_\Iscr$ that is
$c_j=\sprev_l(i)|_\Iscr$ and for any $k\ :\ j<k<i$ one has that $c_k>c_i$ or $c_k<c_j.$
From our assumption one has $\next_r(j)|_\Iscr=\snext_r(j)|_\Iscr=c_i$ and 
$\snext_l(j)|_\Iscr>c_i.$ Thus one has $p(j\pr,\next)|_{\psi(\Iscr)}=i\pr$ and
$ \next_l(j\pr)|_{\psi(\Iscr)}<\snext_r(j\pr)|_{\psi(\Iscr)}.$ 
Thus by (i) $T_{\psi(\Iscr)}(j\pr,\next)\in\Dscr(T_{\psi(\Iscr)}).$ Applying $\psi$ we get
$\psi(T_{\psi(\Iscr)}(j\pr,\next))\in\Dscr(T_\Iscr).$ On the other hand by \ref{3.9} 
$\psi(T_{\psi(\Iscr)}(j\pr,\next))=T_\Iscr(c_j\diar C_i(c_j))=T_\Iscr(C_i,\prev).$
Hence $T_\Iscr(C_i,\prev)\in \Dscr(T_\Iscr).$
\QED
\section
{\bf  On geometry of  descendants of a Richardson orbital variety}
\subsection{}\label{4.1}
As it was shown in [Me] one has $\ov{\gn\cap\Oscr}=\gn\cap\ov\Oscr$
for any $\Oscr.$ The question is whether considering $\gm_{\sr \Iscr}$ instead of
$\gn$ gives the same equality. Of course the first restriction is $\Oscr\geq \Oscr_\Iscr$
since if $\Oscr<\Oscr_\Iscr$ then $\ov{\gm_{\sr\Iscr}\cap\Oscr}=\emptyset$ and
$\gm_{\sr\Iscr}\cap\ov\Oscr=\gm_{\sr \Iscr}$ so the equality  trivially fails to be true.
Let us show that even for $\Oscr>\Oscr_\Iscr$ the equality does not necessarily hold.
Recall that 
$$\gm_{\sr \Iscr}\cap\ov\Oscr=\coprod\limits_{\Oscr\pr\geq \Oscr}\gm_{\sr\Iscr}\cap\Oscr\pr.$$
Now take $\Iscr\st\Pi$ for which there exist $\Vscr\in\Dscr_G(\Vscr_\Iscr)$
such that $\Oscr_\Vscr$ is not a descendant of $\Oscr_\Iscr.$
Let $\Oscr$ be an intermediate orbit, that is $\Oscr_\Iscr<\Oscr<\Oscr_\Vscr.$
Obviously 
$\Vscr\st\gm_{\sr \Iscr}\cap\ov\Oscr.$ On the other hand 
one has by \ref{2.1}
$$\gm_{\sr \Iscr}\cap\Oscr=\cup\Vscr_i,$$
where $\Vscr_i$ are orbital varieties.
Thus $\ov{\gm_{\sr \Iscr}\cap\Oscr}=\cup\ov\Vscr_i.$
Now by \ref{2.6} one has $\Vscr\not\st\ov\Vscr_i$ for any $i$ thus 
$\Vscr\not\st
\cup\ov\Vscr_i$ just by irreducibility of $\Vscr.$
\par
The first example occurs in $\gs\gl_{\sr 4}.$
Let 
$$T_\Iscr=
\vcenter{
\halign{& \hfill#\hfill
\tabskip4pt\cr
\multispan{7}{\hrulefill}\cr
\ssa
\vb & 1 &  & 3 && 4 &\ts\vb\cr
\vsa
&&\multispan{5}{\hrulefill}\cr
\ssa
\vb & 2 & \ts\vb\cr
\vsa
\multispan{3}{\hrulefill}\cr}}\ ,\quad
S=\vcenter{
\halign{& \hfill#\hfill
\tabskip4pt\cr
\multispan{5}{\hrulefill}\cr
\ssa
\vb & 1 &  & 3 &\ts\vb\cr
\vsa
&&&&\cr
\ssa
\vb & 2 && 4 &\ts\vb\cr
\vsa
\multispan{5}{\hrulefill}\cr}}\quad{\rm and}\quad
T=\vcenter{
\halign{& \hfill#\hfill
\tabskip4pt\cr
\multispan{5}{\hrulefill}\cr
\ssa
\vb & 1 &  & 4 &\ts\vb\cr
\vsa
&&\multispan{3}{\hrulefill}\cr
\ssa
\vb & 2 &\ts\vb\cr
\vsa
&&\cr
\ssa
\vb & 3 &\ts\vb\cr
\vsa
\multispan{3}{\hrulefill}\cr}}.
$$
Take $\lam=(2,2).$ One has $\Vscr_T\st\gm_{\sr\Iscr}\cap\ov\Oscr_\lam$
and $\Vscr_T\not\st\ov{\gm_{\sr\Iscr}\cap\Oscr_\lam}=\ov\Vscr_S.$
\par
On the other hand note that by  \ref{2.1} and the definition of 
$\Dscr_G(T_\Iscr)$
one has 
$$\ov{\gm_{\sr \Iscr}\bigcap\left(\bigcup\limits_{T\in\Dscr_G(T_\Iscr)}\Oscr_{\sh T}\right)}=\gm_{\sr \Iscr}\bigcap 
\left(\bigcup\limits_{T\in\Dscr_G(T_\Iscr)}\ov\Oscr_{\sh T}\right)$$
In particular if $|T_\Iscr^1|=2$ or $|T_\Iscr^1|=|T_\Iscr^2|=\ldots=|T_\Iscr^k|$
that is in the cases when all $T\in \Dscr_G(T_\Iscr)$ are of the same shape, call it $\lam,$
one has $\ov{\gm_{\sr\Iscr}\cap\Oscr_\lam}=\gm_{\sr\Iscr}\cap\ov\Oscr_\lam.$ 
\subsection{}\label{4.2}
As we have explained in \ref{1.10} all the descendants 
of codimension 1 of $\Vscr_{\sr \Iscr}$
are complete intersections. 
According to \ref{2.6} some of the descendants of $\Vscr_{\sr \Iscr}$
are of codimension greater than 1 in it. 
The question is whether a descendant of codimension greater than 1
of $\Vscr_{\sr \Iscr}$ is necessarily a complete intersection.
In \ref{4.4} we give an example in $\gs\gl_6$ of $\Vscr_{\sr \Iscr}$
and its descendant of codimension 2 which is not a complete intersection.
\par
First of all note that for $\gs\gl_n$ where $n\leq 5$ all the descendants of
a Richardson orbital variety
are complete intersections. Indeed as it is noted in [vanL, p.16] all the orbital
varieties for $n\leq 4$ are complete intersections and
the only orbital variety in $\gs\gl_{\sr 5}$
which is not a complete intersection is $\Vscr_T$ where
$$T= 
\vcenter{
\halign{& \hfill#\hfill
\tabskip4pt\cr
\multispan{5}{\hrulefill}\cr
\ssa
\vb & 1 &  & 2 &\ts\vb\cr
\vsa
&&&&\cr
\ssa
\vb & 3 &  &  4 & \ts\vb\cr
\vsa
&&\multispan{3}{\hrulefill}\cr
\ssa
\vb & 5 & \ts\vb\cr
\vsa
\multispan{3}{\hrulefill}\cr}}$$
Note that $\tau(T)=\{\al_{\sr 2},\al_{\sr 4}\}=:\Iscr$ and 
that by \ref{2.6} $T\not\in\Dscr_G(T_\Iscr).$
\subsection{}\label{4.3}
Before  considering the example in $\gs\gl_{\sr 6}$ let us introduce the notation we use
in what follows.
\par
Let $e_{i,j}$ be
the matrix having $1$ in the $ij-$th entry
and $0$ elsewhere.
Let $x_{j,i}$ denote the coordinate function on $\gog=\gs\gl_n$ 
defined by
$$x_{j,i}(e_{r,s})=\begin{cases}
-1, & {\rm if}\ (r,s)=(i,j),\cr
 0, & {\rm otherwise}.\cr
\end{cases}$$
Then the Poisson bracket $\{,\}$ defined on 
$S(\gog^*)$ through the Lie bracket on $\gog$ satisfies
$$\{x_{i,j},x_{r,s}\}=\de_{j,r}x_{i,s}-\de_{s,i}x_{r,j} \eqno{(*)} $$
where $\de_{i,j}$ is the Kronecker delta. 
Let $\gm^-$ be the opposed algebra of $\gm.$
Setting $x_{i,j}=e_{i,j}$ 
for $i>j$ identifies $\gm^-$ with $\gm^*.$ Set $A:=S(\gm^-)$ that is a symmetric algebra of $\gm^-.$
\subsection{}\label{4.4}
Consider $\gs\gl_{\sr 6}.$ Set 
$\Iscr=\{\al_{\sr 1},\al_{\sr 2},\al_{\sr 5}\}.$ One has
$$T_\Iscr=
\vcenter{
\halign{& \hfill#\hfill
\tabskip4pt\cr
\multispan{7}{\hrulefill}\cr
\ssa
\vb & 1 &  & 4 && 5&\ts\vb\cr
\vsa
&&&&\multispan{3}{\hrulefill}\cr
\ssa
\vb & 2 &  &  6 & \ts\vb\cr
\vsa
&&\multispan{3}{\hrulefill}\cr
\ssa
\vb & 3 & \ts\vb\cr
\vsa
\multispan{3}{\hrulefill}\cr}}$$
Consider $\Vscr:=\Vscr_T$ where 
$$T=\vcenter{
\halign{& \hfill#\hfill
\tabskip4pt\cr
\multispan{7}{\hrulefill}\cr
\ssa
\vb & 1 &  & 4 && 5&\ts\vb\cr
\vsa
&&\multispan{5}{\hrulefill}\cr
\ssa
\vb & 2 & \ts\vb\cr
\vsa
&&\cr
\ssa
\vb & 3 & \ts\vb\cr
\vsa
&&\cr
\ssa
\vb & 6 & \ts\vb\cr
\vsa
\multispan{3}{\hrulefill}\cr}}$$
Note that $T\in\Dscr_G(T_\Iscr)$ and moreover $\sh(T)=(3,1,1,1)$
is a descendant of $\sh(T_\Iscr)=(3,2,1).$ One has 
$\codim_{\gm_\Iscr}(\Vscr)=2.$
\par
Let $I(\Vscr)$ denote
the ideal of definition of $\Vscr$ in $A.$
Let us show that $I(\Vscr)=\sqrt I$ where $I:=<f_1,f_2,f_3,f_4>$ for
$$\begin{array}{rl}
f_1&=\left\vert\begin{array}{ccc}
x_{\sr 4,2}&x_{\sr 4,3}&0 \cr
x_{\sr 5,2}&x_{\sr 5,3}&x_{\sr 5,4} \cr
x_{\sr 6,2}&x_{\sr 6,3}&x_{\sr 6,4} \cr
\end{array}\right\vert=
x_{\sr 4,2}x_{\sr 5,3}x_{\sr 6,4}+x_{\sr 4,3}x_{\sr 5,4}x_{\sr 6,2}-x_{\sr 4,2}x_{\sr 5,4}x_{\sr 6,3}-
x_{\sr 4,3}x_{\sr 5,2}x_{\sr 6,4}\cr
&\cr
f_2&=\left\vert\begin{array}{ccc}
x_{\sr 4,1}&x_{\sr 4,3}&0 \cr
x_{\sr 5,1}&x_{\sr 5,3}&x_{\sr 5,4} \cr
x_{\sr 6,1}&x_{\sr 6,3}&x_{\sr 6,4} \cr
\end{array}\right\vert=
x_{\sr 4,1}x_{\sr 5,3}x_{\sr 6,4}+x_{\sr 4,3}x_{\sr 5,4}x_{\sr 6,1}-x_{\sr 4,1}x_{\sr 5,4}x_{\sr 6,3}-
x_{\sr 4,3}x_{\sr 5,1}x_{\sr 6,4}\cr
&\cr
f_3&=\left\vert\begin{array}{ccc}
x_{\sr 4,1}&x_{\sr 4,2}&0 \cr
x_{\sr 5,1}&x_{\sr 5,2}&x_{\sr 5,4} \cr
x_{\sr 6,1}&x_{\sr 6,2}&x_{\sr 6,4} \cr
\end{array}\right\vert=
x_{\sr 4,1}x_{\sr 5,2}x_{\sr 6,4}+x_{\sr 4,2}x_{\sr 5,4}x_{\sr 6,1}-x_{\sr 4,1}x_{\sr 5,4}x_{\sr 6,2}-
x_{\sr 4,2}x_{\sr 5,1}x_{\sr 6,4}\cr
&\cr
f_4&=\left\vert\begin{array}{ccc}
x_{\sr 4,1}&x_{\sr 4,2}&x_{\sr 4,3} \cr
x_{\sr 5,1}&x_{\sr 5,2}&x_{\sr 5,3} \cr
x_{\sr 6,1}&x_{\sr 6,2}&x_{\sr 6,3} \cr
\end{array}\right\vert=\cr
&=x_{\sr 4,1}x_{\sr 5,2}x_{\sr 6,3}+x_{\sr 4,2}x_{\sr 5,3}x_{\sr 6,1}+x_{\sr 4,3}x_{\sr 5,1}x_{\sr 6,2}
 -x_{\sr 4,1}x_{\sr 5,3}x_{\sr 6,2}-x_{\sr 4,2}x_{\sr 5,1}x_{\sr 6,3}-x_{\sr 4,3}x_{\sr 5,2}x_{\sr 6,1}\cr
\end{array}$$
Let us show first of all that all four polynomials are irreducible. Indeed $f_1=m_{l_{\Iscr\pr}}$
where $\gm_{\Iscr\pr}=\pi_{2,6}(\gm_\Iscr).$ Thus $f_1$ is irreducible by [JM]. Now $f_2=-\{x_{\sr 2,1},f_1\}$ and 
$f_3=\{x_{\sr 3,1},f_1\}$
so they are irreducible. Finally $f_4=m_{l_\Iscr\prpr}$ where $\Iscr\prpr=\{\al_{\sr 1},\al_{\sr 2},\al_{\sr 4},\al_{\sr 5}\}$
and again is irreducible by [JM]. Further note that $f_1\not\in<f_2,f_3,f_4>$ since 
$f_2|_{x_{4,1}=x_{5,1}=x_{6,1}=0}=f_3|_{x_{4,1}=x_{5,1}=x_{6,1}=0}=f_4|_{x_{4,1}=x_{5,1}=x_{6,1}=0}
= 0$ and $f_1|_{x_{4,1}=x_{5,1}=x_{6,1}=0}\ne 0.$
Exactly in the same way we show that $f_2\not\in<f_1,f_3,f_4>$, $f_3\not\in<f_1,f_2,f_4>$
and $f_4\not\in<f_1,f_2,f_3>.$ Thus $I$ is generated by 4 polynomials. One can see at once that $I$
is stable under Poison bracket action.
\par
Let $\Vscr(I)$ denote variety of $I.$ It consists of
all the matrices in
$\gm_{\sr \Iscr}$ of rank less or equal to 2 that is 
$\Vscr(I)=\gm\cap(\Oscr_{\lam_1}\cup\Oscr_{\lam_2}\cup\Oscr_{\lam_3}\cup\Oscr_{\lam_4})$
where $\lam_1=(3,1,1,1),\ \lam_2=(2,2,1,1),\ \lam_3=(2,1,1,1,1),\ 
\lam_4=(1,1,1,1,1,1).$ In particular $\ov\Vscr_T\st\Vscr(I).$
By \ref{2.1} one has  
$$\begin{array}{rl}
\gm_{\sr \Iscr}\cap\Oscr_{\lam_1}&=\Vscr_T\cr
           \gm_{\sr \Iscr}\cap\Oscr_{\lam_2}&=\Vscr_P\cup\Vscr_Q\cr
           \gm_{\sr \Iscr}\cap\Oscr_{\lam_3}&=\Vscr_S\cup\Vscr_U\cr
           \gm_{\sr \Iscr}\cap\Oscr_{\lam_4}&=\Vscr_Y\cr
\end{array}$$
where
$$P=\vcenter{
\halign{& \hfill#\hfill
\tabskip4pt\cr
\multispan{5}{\hrulefill}\cr
\ssa
\vb & 1 &  & 4 &\ts\vb\cr
\vsa
&&&&\cr
\ssa
\vb & 2 && 5 & \ts\vb\cr
\vsa
&&\multispan{3}{\hrulefill}\cr
\ssa
\vb & 3 & \ts\vb\cr
\vsa
&&\cr
\ssa
\vb & 6 & \ts\vb\cr
\vsa
\multispan{3}{\hrulefill}\cr}},\
Q=\vcenter{
\halign{& \hfill#\hfill
\tabskip4pt\cr
\multispan{5}{\hrulefill}\cr
\ssa
\vb & 1 &  & 5 &\ts\vb\cr
\vsa
&&&&\cr
\ssa
\vb & 2 && 6 & \ts\vb\cr
\vsa
&&\multispan{3}{\hrulefill}\cr
\ssa
\vb & 3 & \ts\vb\cr
\vsa
&&\cr
\ssa
\vb & 4 & \ts\vb\cr
\vsa
\multispan{3}{\hrulefill}\cr}},\
S=
\vcenter{
\halign{& \hfill#\hfill
\tabskip4pt\cr
\multispan{5}{\hrulefill}\cr
\ssa
\vb & 1 &  & 4 &\ts\vb\cr
\vsa
&&\multispan{3}{\hrulefill}\cr
\ssa
\vb & 2 & \ts\vb\cr
\vsa
&&\cr
\ssa
\vb & 3 & \ts\vb\cr
\vsa
&&\cr
\ssa
\vb & 5 & \ts\vb\cr
\vsa
&&\cr
\ssa
\vb & 6 & \ts\vb\cr
\vsa
\multispan{3}{\hrulefill}\cr}},\
U=\vcenter{
\halign{& \hfill#\hfill
\tabskip4pt\cr
\multispan{5}{\hrulefill}\cr
\ssa
\vb & 1 &  & 5 &\ts\vb\cr
\vsa
&&\multispan{3}{\hrulefill}\cr
\ssa
\vb & 2 & \ts\vb\cr
\vsa
&&\cr
\ssa
\vb & 3 & \ts\vb\cr
\vsa
&&\cr
\ssa
\vb & 4 & \ts\vb\cr
\vsa
&&\cr
\ssa
\vb & 6 & \ts\vb\cr
\vsa
\multispan{3}{\hrulefill}\cr}}\ {\rm and}\quad
Y=\vcenter{
\halign{& \hfill#\hfill
\tabskip4pt\cr
\multispan{3}{\hrulefill}\cr
\ssa
\vb & 1 &\ts\vb\cr
\vsa
&&\cr
\ssa
\vb & 2 & \ts\vb\cr
\vsa
&&\cr
\ssa
\vb & 3 & \ts\vb\cr
\vsa
&&\cr
\ssa
\vb & 4 & \ts\vb\cr
\vsa
&&\cr
\ssa
\vb & 5 & \ts\vb\cr
\vsa
&&\cr
\ssa
\vb & 6 & \ts\vb\cr
\vsa
\multispan{3}{\hrulefill}\cr}}.$$
By Part I, 3.3.3 one can see that $P,Q,S,U,Y\dgs T.$ Thus $\Vscr(I)\st \ov\Vscr_T.$
We get $\Vscr(I)= \ov\Vscr_T.$  
\par
To show that $\ov\Vscr_T$ is not a complete intersection we
must show that $\sqrt I$ is generated by at least 3
polynomials. We will show that it is generated by at least 4 polynomials.
\par
Let $\sqrt I=<p_i>_{i=1}^k.$ Let us write  $p_i=\sum p_i^j$ where $p_i^j$ is a component
of $p_i$ whose monomials contain $j$ different variables.
\par 
First let us show that $p_i^j= 0$ for $j\leq 2.$ 
\begin{itemize}
\item[(i)] It is obvious that $p_i^{\sr 0}= 0$
since there exists $t_i$ such that $p_i^{t_i}\in I.$  
\item[(ii)] Note that $A$ has 11 variables $x_{i,j}.$
Note that for any specialization $\ov x$ of
$\{x_{i,j}\}$ taking 10 variables to 0 one has $f_i|_{\ov x}=0$ for $i=1,2,3,4.$
Thus $I|_{\ov x}=0.$
On the other hand $p_i|_{\ov x}=p_i^1|_{\ov x}$ and $p_i^{t_i}|_{\ov x}=(p_i^1)^t_i|_{\ov x}\in I|_{\ov x}.$
Thus $p_i^1|_{\ov x}=0$ for any such specialization which provides $p_i^1= 0.$
\item[(iii)] Now note that for any specialization $\ov x$ of
$\{x_{i,j}\}$ taking 9 variables to 0 one has $f_i|_{\ov x}=0$ for $i=1,2,3,4.$
Thus exactly in the same way as in (ii) we get that $p_i^2=0.$
\end{itemize}
Now let $p_{i,j}$ denote the component of $p_i$ of total degree $j.$ 
Then $p_i=p_{i,\sr 3}+\hat p_i$ where $\hat p_i=\sum_{j\geq 4}p_{i,j}.$ 
Since $f_1,f_2,f_3,f_4\in \sqrt I$ we get by degree considerations
$$f_j=\sum_i B_i^j p_i=\sum_i b_i^j p_{i,3}$$
where $B_i^j$ are some polynomials and $b_i^j\in\Co.$
Consider a linear affine space $V={\Co}^r$
with the basis $\{x_{i,j}x_{p,q}x_{s,t}\}$ where $(i,j),\ (p,q),\ (s,t)$
are pairwise different and $r=\binom{11}{3}.$
Since $f_1,\ldots,f_4$ are linearly independent in $V$
one has that the number of $p_i$ must be at least 4.
\par
Moreover using more subtle analysis one can show that for any 
$i\ :\ 1\leq i\leq k$ and any $j\ :\ j\geq 3$ one has
$p_{i,j}\in\sqrt I$ and then using  specializations
taking $\sqrt I$ onto $I_i=<f_i>$ which are all simple ideals
and some further specializations
one can show that $\sqrt I=<f_1,f_2,f_3,f_4>.$
\bigskip
\parno
\centerline{ INDEX OF NOTATION}

\parno
\begin{tabular}{llll}
\ref{1.2} & $\Oscr,\ \Vscr,\ \Oscr_{\Vscr}$&\ref{2.01}
          & $T_j^i,\ \omega^i(T),\ T_j,\ T^i,\ T^{i,j},\ T^{i,\infty},\ |(T_j)|,$\cr
\ref{1.3} & $R,\ R^+,\ \Pi,\ X_{\al},\ \bB,\ \gb,\  \Vscr_w,\ \Oscr_w$&
          & $|(T^i)|,\ <T>,\ {\boT}_n,\ {\boT}_E,\ (T,S),\ T^{\dagger}$\cr
\ref{1.4} & $\Iscr,\ {\bP}_{\Iscr},\ {\bM}_{\Iscr},\ \gm_{\sr \Iscr},$&\ref{2.2} 
          & $T_{\Iscr},\ C^{\Iscr}_i,\ c_i,\ \sigma_i$\cr
          & $W_{\Iscr},\ w_{\sr \Iscr},\ \Vscr_{\Iscr},\ \Oscr_{\Iscr}$&\ref{2.4} 
          & $\hat T_{s,t},\ C^{1,j-1}_i,\ C_i(j),\ T_{\Iscr}(c_s\diar C_i(j))$\cr
\ref{1.5} & $\lambda,\ \Oscr_{\lambda},\ \lambda\leq \mu,\ \Oscr\leq \Oscr\pr$&\ref{2.5} 
          & $\next_l(i),\ \prev_l(i),\ \sprev_l(i),\ \snext_r(i),$\cr
\ref{1.6} & $\Cscr_{\Vscr},\ \Cscr_w,\ \bS_n,\  T(w),$&
          & $\sprev_r(i),\ T_{\Iscr}(i,\next),\ T_{\Iscr}(i,\prev), T_{\Iscr}(i)$\cr
          & $\bT_n,\ \Vscr_T,\ T_{\Vscr},\ \sh(T)$ & \ref{2.6} & $\Dscr_G(T)$\cr         
\ref{1.6a}& $r_{\sr T}(u),\ \tau(T),\ \tau(\Vscr)$ & \ref{3.0} & $<i,j>,\ \pi_{i,j}$\cr
\ref{1.7} & $\Vscr\go \Wscr$&\ref{3.1} 
          & $\Iscr_n,\ \Iscr_1,\ p(i,\next),\ p(i,\prev)$\cr
\ref{2.0} & $S(w),\ y\dor w,\ \tau(w),\ {\bP}_{\al},\ \al_{i,j},$&\ref{3.3} & ${\binom{P}{Q}}$\cr
          & $p_w(i),\ <w>,\ \ov w,\ [x,y],\ {\boS}_n,\ {\boS}_E$&\ref{4.3} & $x_{i,j},\ \gm^-,\ A$\cr
\end{tabular}
\bigskip
\parno

\end{document}